\documentclass[pdflatex,sn-mathphys-num]{sn-jnl}%

\usepackage{graphicx}%
\usepackage{multirow}%
\usepackage{amsmath,amssymb,amsfonts}%
\usepackage{amsthm}%
\usepackage{mathrsfs}%
\usepackage[title]{appendix}%
\usepackage{xcolor}%
\usepackage{booktabs}%
\usepackage{cleveref}%
\usepackage{array}
\crefname{equation}{}{}

\begin{document}

\title[Failure Modes of IPM SDP Solvers on Moment--SOS Hierarchies]{Failure Modes of Interior-Point SDP Solvers on Moment--SOS Hierarchies}

\author*[1,2]{\fnm{Michal} \sur{Ko\v{c}vara}}\email{m.kocvara@bham.ac.uk}

\affil[1]{\orgdiv{School of Mathematics}, \orgname{University of Birmingham}, \orgaddress{\city{Birmingham}, \postcode{B15 2TT}, \country{United Kingdom}}}

\affil[2]{\orgdiv{Faculty of Civil Engineering}, \orgname{Czech Technical University in Prague}, \orgaddress{\street{Th\'akurova 2077/7}, \city{Prague~6}, \postcode{16629}, \country{Czech Republic}}}

 \abstract{\unboldmath %
	Standard semidefinite programming (SDP) solvers working in double precision can become
	unreliable on difficult problems, even when their output appears plausible. We
	study this issue for SDP problems arising from high-order moment-SOS
	hierarchies. Three polynomial optimization problems are analysed by comparing
	standard double-precision computations with reference solutions obtained
	using the extended-precision solver \mbox{Loraine\,64x4}. We observe two
	distinct patterns of numerical failure. In \emph{false convergence}, the
	computed bounds form a smooth sequence but approach the limit at the wrong
	rate. In \emph{stagnation}, the sequence develops a plateau or oscillations
	and no longer reflects the convergence of the underlying hierarchy. We
	introduce an empirical effective condition number based on the positive
	spectral components of the primal and dual solution matrices. It serves as a
	warning indicator rather than an error bound, and its growth accompanies
	numerical deterioration in our examples. We also examine the imbalance
	between primal and dual spectral conditioning and a Cholesky inverse residuum
	as potential diagnostic indicators. The results show that solver termination
	messages and apparently regular convergence curves may be insufficient to
	establish the reliability of high-order moment-SOS bounds.}



\maketitle

\section{Introduction}\label{sec:intro}

The moment-SOS hierarchy generates a sequence of semidefinite programming
(SDP) relaxations whose optimal values converge to the global minimum of
a polynomial optimization problem (POP) as the relaxation order $\omega$
increases; see, e.g., \cite{lasserre2009moments}. Although some polynomial
optimization problems admit finite convergence, the three problems studied
here exhibit only asymptotic convergence. In finite precision, the
higher-order relaxations can become increasingly ill-conditioned, and
their numerical solution can lose accuracy.

In this paper, we study two qualitatively different manifestations of
this loss of accuracy. The first is \emph{false convergence}: the computed
relaxation values follow a smooth and apparently convincing convergence
pattern, but with an incorrect asymptotic rate. The second is
\emph{stagnation}: the computed values cease to follow a systematic
convergence trend and instead exhibit a plateau or oscillations.

We analyse three polynomial optimization problems exhibiting these
behaviours. Problems~1 and~2 show false convergence, while Problem~3 shows
stagnation. Our analysis focuses on two distinct aspects of numerical
accuracy. The first is the spectral conditioning of the primal and dual
solution matrices, which affects the accuracy of the Newton step in an
interior-point (IP) method. The second is the amplification of objective
error when the optimal value is small. We quantify the latter effect by
an objective amplification factor $C$, defined in \Cref{sec:objective}. In our
examples, this factor grows by several orders of magnitude with the
relaxation order, so that even a modest error in the computed objective
relative to the scale of the underlying matrices can result in a large
relative error in the reported bound.

To investigate the spectral conditioning, we analyse the Newton step of
a primal-dual IP method, focusing on the formation and solution of the
Schur complement system. This motivates the definition of the effective condition number
\[
\kappa_{\mathrm{eff}}=\kappa(X_+)\kappa(S_+),
\]
where $X_+$ and $S_+$ denote the restrictions of the exact primal and
dual solution matrices to their positive eigenspaces. The product
provides a heuristic measure of the combined precision demand associated
with these spectral components. Its growth is associated with the
observed deterioration in our examples, but, due to the limited number of problems studied, it should not be considered a general condition number. In
numerical experiments, these condition numbers must be estimated from
sufficiently accurate computed solutions.

We also consider the ratio
\[
r=\frac{\kappa(S_+)}{\kappa(X_+)},
\]
which measures the imbalance between the dual and primal spectral
conditioning. In our three examples, the false-convergence cases are
strongly dual-dominated, whereas the stagnation case has comparable
primal and dual factors. Of course, three examples are not enough to
establish a general classification, but the pattern is worth noting.

Computing these diagnostics reliably may require a high-precision solution, both to obtain accurate spectral information and to establish a trustworthy reference objective. These may not always be available or practical.
We therefore also investigate a cheaper diagnostic that can be evaluated
during the IP iteration: the Cholesky inverse residuum
\[
\rho=\|SS^{-1}_{\mathrm{comp}}-I\|,
\]
where $S^{-1}_{\mathrm{comp}}$ denotes the inverse action on $S$ computed using
triangular solves, without forming an explicit inverse. In solvers that
support multiple arithmetic precisions, monitoring $\rho$ may provide an
early warning that higher precision is needed. The aim is to switch
precision while the Schur complement matrix is still positive definite,
rather than waiting for it to become indefinite. The appropriate
threshold is solver- and problem-dependent, and $\rho$ does not monitor
every source of numerical error.

Numerical difficulties of SDP solvers on moment-SOS relaxations have
been reported previously. Waki, Nakata and Muramatsu~\cite{waki2012}
documented incorrect bounds, non-convergence, and other anomalous behaviour
of interior-point methods applied to polynomial optimization, and studied
the effect of extended-precision SDP solvers. The present paper sheds light on the numerical mechanisms behind false convergence and stagnation in our examples.

\section{Numerical observations}\label{sec:observations}

\subsection{Polynomial optimization and moment-SOS relaxations}

We first recall a primal-dual SDP pair in the sign convention used below:
\begin{equation}\label{eq:sdp-pair}
\begin{alignedat}{2}
\text{primal:}\quad
&\max_{X\in\mathbb S^m} && \langle B,X\rangle\\
&\text{subject to}\ && \mathcal A^*(X)-c=0,\\
&&& X\succeq 0,
\end{alignedat}
\qquad
\begin{alignedat}{2}
\text{dual:}\quad
&\min_{y\in\mathbb R^n,\,S\in\mathbb S^m}
&& c^Ty\\
&\text{subject to}\ 
&& \mathcal A(y)-B-S=0,\\
&&& S\succeq 0.
\end{alignedat}
\end{equation}
Here, for any $U,V\in \mathbb S^m$, $\langle U,V\rangle=\operatorname{tr}(UV)$ is the Frobenius inner
product, $\mathcal A(y)=\sum_i y_iA_i$, and
$\mathcal A^*(X)_i=\langle A_i,X\rangle$ with $A_i\in\mathbb S^m$, $i=1,\ldots,n$. Note that in the actual moment-SOS
formulations there can be several positive semidefinite constraints which, formally, can be written as \cref{eq:sdp-pair} with block diagonal matrices $X$ and $S$ of the same block structure. With some abuse of notation, unless stated otherwise, in the text below
$X$ and $S$ will only denote the principal Gram and moment blocks for
which spectral data will be retained.

%

The objective of polynomial optimization (POP) is to determine the optimal value
\[
f^* := \min_{x\in\mathbb{R}^\nu} f_0(x) \quad\text{subject to}\quad f_i(x)\ge 0,\ i=1,\ldots,\mu,
\]
where $f_0,\ldots,f_\mu$ are real polynomials. The moment-SOS hierarchy
approximates~$f^*$ from below by a sequence of SDPs indexed by the
relaxation order~$\omega$. The dual variable~$S$ is a truncated moment
matrix, together with a localizing matrix for each constraint~$f_i$, while
the primal variable~$X$ consists of Gram matrices of an
SOS certificate 
$$
f_0-f_\omega = \sigma_0+\sum_{i=1}^\mu\sigma_i f_i,
$$
where $\sigma_0,\ldots,\sigma_\mu$ are SOS polynomials. The optimal
value $f_\omega$ of the relaxation satisfies $f_\omega\le f^*$ for every~$\omega$ and, under standard assumptions, converges monotonically to $f^*$ as $\omega\to\infty$. For details, see \cite{lasserre2015,Parrilo:2000thesis,Nie2023}.

The problems studied below share two key features: 
\begin{itemize}
	\item[--] their optimal objective value is zero, and
	\item[--] the moment-SOS hierarchy does not attain this value at any finite relaxation order~$\omega$.
\end{itemize} 
As~$\omega$ increases, the resulting sequence of SDPs grows in dimension, and its optimal values converge monotonically to zero without ever reaching it. 
Notice that the moment-SOS finite-convergence result requires standard
nonlinear-programming conditions, LICQ, strict complementarity, and
second-order sufficient conditions to be satisfied; see \cite[Thm.~6.5]{lasserre2015}. Some of these conditions are violated in each of our problems.

The examples below are taken from \cite{baldi-kocvara}, where our primary interest was the rate of convergence. Numerical experiments aimed at estimating this rate revealed the stagnation and false convergence described above. In this article, we investigate these phenomena by comparing standard double-precision SDP solvers with the solver Loraine running in high-precision arithmetic.

\subsection{Solvers and reference computations}\label{sec:solvers}
To solve the resulting SDP problems, we use two popular SDP solvers, MOSEK~\cite{mosek} and CSDP~\cite{csdp}  called through Julia/JuMP with their default
settings. We compare their results with that of the high-precision Julia package {\tt Loraine.jl}~\cite{loraine}\footnote{\url{https://github.com/kocvara/Loraine.jl}}. Loraine relies on the extended-precision arithmetic provided by the Julia package {\tt MultiFloats.jl}. In addition to the standard {\tt Float64} type, it defines types {\tt Float64x[N]} with ${\tt N}=2,\ldots,8$. For the purpose of this study, the values computed by Loraine running in Float64x4 extended precision are treated as the reference solution. Below we refer to this solver as \mbox{Loraine\,64x4}.

The Loraine
stopping criterion is based on the DIMACS error measures \cite{dimacs}. In the
sign convention of \eqref{eq:sdp-pair}, these are
\begin{align*}
 \mathrm{err}_1&=\frac{\|\mathcal A^*(X)-c\|}{1+\|c\|},&
 \mathrm{err}_2&=\frac{\max\{0,-\lambda_{\min}(X)\}}{1+\|c\|},\\
 \mathrm{err}_3&=\frac{\|\mathcal A(y)-B-S\|}{1+\|B\|},&
 \mathrm{err}_4&=\frac{\max\{0,-\lambda_{\min}(S)\}}{1+\|B\|},\\
 \mathrm{err}_5&=\frac{c^Ty-\langle B,X\rangle}
 {1+|\langle B,X\rangle|+|c^Ty|},&
 \mathrm{err}_6&=\frac{\langle S,X\rangle}
 {1+|\langle B,X\rangle|+|c^Ty|}.
\end{align*}
In the numerical experiments below, \mbox{Loraine\,64x4} terminates when
$$
\varepsilon_{_{\rm DIMACS}} := \max\{\mbox{err}_1 ,\mbox{err}_2 ,\mbox{err}_3 ,\mbox{err}_4 ,\mbox{err}_5,\mbox{err}_6\}\leq 10^{-16}. 
$$
Some spectra are also compared with runs stopped at $10^{-6}$ to show that
a small aggregate residuum need not imply reliable numerical rank.

The SDP problems are generated by
{\tt SumOfSquares.jl} and then dualized before being passed to Loraine via
\[
\verb!dual_loraine = Dualization.dual_optimizer(loraine)!.
\]
This dualization is required for the efficiency of Loraine. In this setting, the moment matrix is the dual variable $S$ in
\cref{eq:sdp-pair} and the Gram matrix is the primal variable~$X$.

\subsection{Problem 1, false convergence}\label{sec:ex2a}
We first consider the univariate polynomial optimization problem due to Stengle \cite{stengle},
\begin{equation}\label{eq:p1}
\min_{x} 1-x^2\quad \text{subject to}\quad (1-x^2)^3 \geq 0
\end{equation}
and solve it by the moment-SOS hierarchy. 
Because of the formulation of the constraint, LICQ fails at the optimal points  $x^*=\pm 1$. Therefore, finite convergence is not guaranteed by the standard theory and, in fact, the hierarchy converges only asymptotically as shown in \cite{stengle}.
Moreover, for this problem the optimal values of the SDP relaxations are known in closed form: Henrion \cite{didier} proved that, for relaxation order~$\omega$, the optimal value of the SDP relaxation is 
\begin{equation}\label{eq:didi}
f_\omega = -\frac{1}{\omega(\omega-2)}.
\end{equation}

This problem exhibits the more dangerous failure mode: rather than stagnating for increasing $\omega$, the relaxations computed by the standard solvers appear to converge at a convincing rate that is, however, entirely spurious.

In \Cref{tab:prob1obj}, we  compare \mbox{Loraine\,64x4} with CSDP and MOSEK using double-precision arithmetic. The results of \mbox{Loraine\,64x4} agree with the reference rate of \cref{eq:didi} to all digits reported in the table. Because the CSDP primal and dual objectives 
noticeably disagree for larger~$\omega$, both are listed in the table. For MOSEK we only show primal objectives as a confirmation of the trend.
{\small
	\begin{table}[h]
		\centering
				\caption{Problem 1, objective values computed by \mbox{Loraine\,64x4}, CSDP and MOSEK}
		\label{tab:prob1obj}
		\begin{tabular}{rlccc}
			\hline
			$\omega$&
			Loraine\,64x4&	CSDP primal & CSDP dual &	MOSEK\\
			\hline
			5 & $-6.66666666666665 \!\!\times\!\! 10^{-2}$	& $-6.66668 \!\!\times\!\! 10^{-2}$ & $-6.66668 \!\!\times\!\! 10^{-2}$ & $-6.66666 \!\!\times\!\! 10^{-2}$  \\
			6 & $-4.16666666666667 \!\!\times\!\! 10^{-2}$ & $-4.16666 \!\!\times\!\! 10^{-2}$ & $-4.16673 \!\!\times\!\! 10^{-2}$ & $-4.16666 \!\!\times\!\! 10^{-2}$    \\
			7 & $-2.85714285714285 \!\!\times\!\! 10^{-2}$	& $-2.85715 \!\!\times\!\! 10^{-2}$ & $-2.86242 \!\!\times\!\! 10^{-2}$ & $-2.85714 \!\!\times\!\! 10^{-2}$  \\
			8 & $-2.08333333333333 \!\!\times\!\! 10^{-2}$ & $-2.08314 \!\!\times\!\! 10^{-2}$ & $-2.12137 \!\!\times\!\! 10^{-2}$ & $-2.22536 \!\!\times\!\! 10^{-2}$  \\
			9 & $-1.58730158730158 \!\!\times\!\! 10^{-2}$	& $-1.64439 \!\!\times\!\! 10^{-2}$ & $-1.80615 \!\!\times\!\! 10^{-2}$ & $-2.06113 \!\!\times\!\! 10^{-2}$  \\
			10 & $-1.25000000000000 \!\!\times\!\! 10^{-2}$ & $-1.70738 \!\!\times\!\! 10^{-2}$ & $-1.82583 \!\!\times\!\! 10^{-2}$ & $-1.64743 \!\!\times\!\! 10^{-2}$  \\
			11 & $-1.01010101010099 \!\!\times\!\! 10^{-2}$	& $-1.50450 \!\!\times\!\! 10^{-2}$ & $-1.65009 \!\!\times\!\! 10^{-2}$ & $-1.52510 \!\!\times\!\! 10^{-2}$  \\
			13 & $-6.99300699300696 \!\!\times\!\! 10^{-3}$	& $-1.18870 \!\!\times\!\! 10^{-2}$ & $-1.30143 \!\!\times\!\! 10^{-2}$ & $-1.28770 \!\!\times\!\! 10^{-2}$  \\
			15 & $-5.12820512820488 \!\!\times\!\! 10^{-3}$	& $-9.88454 \!\!\times\!\! 10^{-3}$ & $-1.06026 \!\!\times\!\! 10^{-2}$ & $-1.07471 \!\!\times\!\! 10^{-2}$  \\
			17 & $-3.92156862745087 \!\!\times\!\! 10^{-3}$	& $-9.74368 \!\!\times\!\! 10^{-3}$ & $-1.07927 \!\!\times\!\! 10^{-2}$ & $-9.59489 \!\!\times\!\! 10^{-3}$  \\
			19 & $-3.09597523219807 \!\!\times\!\! 10^{-3}$	& $-7.58612 \!\!\times\!\! 10^{-3}$ & $-8.20946 \!\!\times\!\! 10^{-3}$ & $-8.47775 \!\!\times\!\! 10^{-3}$  \\
			21 & $-2.50626566416014 \!\!\times\!\! 10^{-3}$	& $-6.58259 \!\!\times\!\! 10^{-3}$ & $-7.13888 \!\!\times\!\! 10^{-3}$ & $-7.84202 \!\!\times\!\! 10^{-3}$ \\
			23 & $-2.07039337474112 \!\!\times\!\! 10^{-3}$	& $-6.38502 \!\!\times\!\! 10^{-3}$ & $-6.94889 \!\!\times\!\! 10^{-3}$ & $-6.85201 \!\!\times\!\! 10^{-3}$ \\
			25 & $-1.73913043478250 \!\!\times\!\! 10^{-3}$	& $-5.61803 \!\!\times\!\! 10^{-3}$ & $-6.11551 \!\!\times\!\! 10^{-3}$ & $-6.55448 \!\!\times\!\! 10^{-3}$ \\
			27 & $-1.48148148148142 \!\!\times\!\! 10^{-3}$	& $-5.82504 \!\!\times\!\! 10^{-3}$ & $-6.39292 \!\!\times\!\! 10^{-3}$ & $-5.90841 \!\!\times\!\! 10^{-3}$ \\
			29 & $-1.27713920817352 \!\!\times\!\! 10^{-3}$	& $-5.05265 \!\!\times\!\! 10^{-3}$ & $-5.52473 \!\!\times\!\! 10^{-3}$ & $-5.73612 \!\!\times\!\! 10^{-3}$ \\
			31 & $-1.11234705227999 \!\!\times\!\! 10^{-3}$	& $-4.94314 \!\!\times\!\! 10^{-3}$ & $-5.42774 \!\!\times\!\! 10^{-3}$ & $-5.92840 \!\!\times\!\! 10^{-3}$ \\
			33 & $-9.77517106549363 \!\!\times\!\! 10^{-4}$	& $-4.63170 \!\!\times\!\! 10^{-3}$ & $-5.08971 \!\!\times\!\! 10^{-3}$ & $-5.28825 \!\!\times\!\! 10^{-3}$ \\
			35 & $-8.65800865800866 \!\!\times\!\! 10^{-4}$	& $-4.32800 \!\!\times\!\! 10^{-3}$ & $-4.75668 \!\!\times\!\! 10^{-3}$ & $-4.81734 \!\!\times\!\! 10^{-3}$ \\
			37 & $-7.72200772200772 \!\!\times\!\! 10^{-4}$	& $-4.25650 \!\!\times\!\! 10^{-3}$ & $-4.69989 \!\!\times\!\! 10^{-3}$ & $-4.45704 \!\!\times\!\! 10^{-3}$ \\
			39 & $-6.93000693000693 \!\!\times\!\! 10^{-4}$	& $-4.11928 \!\!\times\!\! 10^{-3}$ & $-4.56013 \!\!\times\!\! 10^{-3}$ & $-4.57947 \!\!\times\!\! 10^{-3}$ \\
			41 & $-6.25390869293308 \!\!\times\!\! 10^{-4}$	& $-3.83662 \!\!\times\!\! 10^{-3}$ & $-4.23706 \!\!\times\!\! 10^{-3}$ & $-4.25095 \!\!\times\!\! 10^{-3}$  \\
			43 & $-5.67214974475326 \!\!\times\!\! 10^{-4}$	& $-3.65828 \!\!\times\!\! 10^{-3}$ & $-4.04577 \!\!\times\!\! 10^{-3}$ & $-4.09481 \!\!\times\!\! 10^{-3}$  \\
			45 & $-5.16795865633075 \!\!\times\!\! 10^{-4}$	& $-3.89736 \!\!\times\!\! 10^{-3}$ & $-4.37936 \!\!\times\!\! 10^{-3}$ & $-4.00572 \!\!\times\!\! 10^{-3}$  \\
			\hline                               
		\end{tabular} 
	\end{table}
}%
\begin{figure}[h]
	\centering
	{\includegraphics[width=0.6\linewidth]{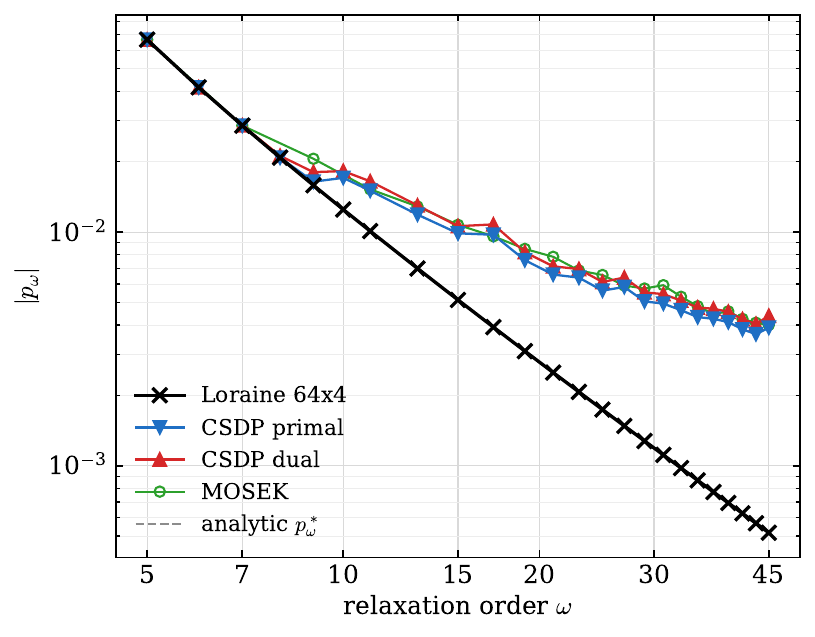}}
	\caption{Problem~1: computed objective $|f_\omega|$ versus relaxation order
		(log-log).  Series as in the legend.}
	\label{fig:prob1}
\end{figure}

\Cref{fig:prob1} shows $|f_\omega|$ versus $\omega$.  At low orders the
three solvers agree.  From about $\omega=9$ onward, CSDP and MOSEK continue to
produce smooth decreasing sequences, although their termination messages no
longer indicate successful optimality.  For the CSDP dual objective, a
least-squares fit over $\omega\in[9,45]$ has slope $-0.98$ ($R^2=0.98$),
whereas the exact rate from \cref{eq:didi} is $-2$.  The smooth but incorrect
rate is the false-convergence behaviour.

The final messages of CSDP, MOSEK and other SDP solvers installed in Julia/JuMP for $\omega=11$ are given in Appendix~\ref{app:solvers}. They show that some solvers even claim optimality despite returning an incorrect optimal objective value.

\subsection{Problem 2, false convergence}\label{sec:ex3}
Next, we consider the polynomial optimization problem in~$\mathbb{R}^2$:
\begin{equation}\label{eq:p3}
\min_{x\in{\mathbb R}^2}\ (1-x_1^2)(1-x_2^2) \quad \text{subject to} \quad
1-x_1^2\ge 0,\quad 1-x_2^2\ge 0\,.
\end{equation}
The feasible set is the square $[-1,1]^2$, and the true optimal value
$f^*=0$ is attained on the entire boundary
$\{(x\in [-1,1]^2\mid |x_1|=1\text{ or } |x_2|=1\}$. At the corners, where both constraints are active, the corresponding multipliers are equal to zero, hence strict complementarity of the nonlinear program does not hold and
so finite convergence is not guaranteed. In~\cite{baldi-kocvara}, it is conjectured that the exact optimal values of the
relaxations are given by
\begin{equation}\label{eq:p3formula}
f_\omega = -\frac{\alpha^4/4}{(\omega-c)^2\bigl[(\omega-c)^2-\alpha^2\bigr]},
\quad
c=\begin{cases}1, & \omega\ \text{odd}\\[2pt]
\alpha, & \omega\ \text{even}\end{cases},\quad
\alpha = \frac{2}{1 + \sqrt{2}}\,.
\end{equation}
This conjecture has been confirmed numerically using high-precision
arithmetic. We use formula \cref{eq:p3formula} as a reference in the comparisons below.

\Cref{tab:prob2obj} compares the objective values computed by \mbox{Loraine\,64x4} with those
obtained by CSDP and MOSEK.  As in Problem~1, we list both primal and dual objectives of CSDP and only primal objectives for MOSEK, as a confirmation of the trend.
{\small
	\begin{table}[h]
		\centering
				\caption{Problem~2: objective values computed by \mbox{Loraine\,64x4}, CSDP
			(primal and dual), and MOSEK (primal). }
		\label{tab:prob2obj}
		\begin{tabular}{rlccc}
			\hline
			$\omega$ & \multicolumn{1}{c}{\mbox{Loraine\,64x4}} & CSDP primal & CSDP dual & MOSEK\\
			\hline
			3  & $-8.88347648318441 \!\!\times\!\! 10^{-3}$ & $-8.88347 \!\!\times\!\! 10^{-3}$ & $-8.88347 \!\!\times\!\! 10^{-3}$ & $-8.88348 \!\!\times\!\! 10^{-3}$\\
			4  & $-1.56639974087585 \!\!\times\!\! 10^{-3}$ & $-1.56640 \!\!\times\!\! 10^{-3}$ & $-1.56641 \!\!\times\!\! 10^{-3}$ & $-1.56641 \!\!\times\!\! 10^{-3}$\\
			5  & $-4.80570260378324 \!\!\times\!\! 10^{-4}$ & $-4.80570 \!\!\times\!\! 10^{-4}$ & $-4.80590 \!\!\times\!\! 10^{-4}$ & $-4.80652 \!\!\times\!\! 10^{-4}$\\
			6  & $-1.72449392209694 \!\!\times\!\! 10^{-4}$ & $-1.72248 \!\!\times\!\! 10^{-4}$ & $-1.72458 \!\!\times\!\! 10^{-4}$ & $-1.72527 \!\!\times\!\! 10^{-4}$\\
			7  & $-9.26214171150048 \!\!\times\!\! 10^{-5}$ & $-9.10503 \!\!\times\!\! 10^{-5}$ & $-9.40041 \!\!\times\!\! 10^{-5}$ & $-1.12136 \!\!\times\!\! 10^{-4}$\\
			8  & $-4.51163274991384 \!\!\times\!\! 10^{-5}$ & $-4.41290 \!\!\times\!\! 10^{-5}$ & $-5.24373 \!\!\times\!\! 10^{-5}$ & $-5.83430 \!\!\times\!\! 10^{-5}$\\
			9  & $-2.90589236327582 \!\!\times\!\! 10^{-5}$ & $-3.17777 \!\!\times\!\! 10^{-5}$ & $-3.61863 \!\!\times\!\! 10^{-5}$ & $-4.44459 \!\!\times\!\! 10^{-5}$\\
			10 & $-1.67779840418985 \!\!\times\!\! 10^{-5}$ & $-2.14343 \!\!\times\!\! 10^{-5}$ & $-2.23264 \!\!\times\!\! 10^{-5}$ & $-3.14480 \!\!\times\!\! 10^{-5}$\\
			11 & $-1.18562691768419 \!\!\times\!\! 10^{-5}$ & $-2.01324 \!\!\times\!\! 10^{-5}$ & $-1.90779 \!\!\times\!\! 10^{-5}$ & $-2.49631 \!\!\times\!\! 10^{-5}$\\
			12 & $-7.60141972636513 \!\!\times\!\! 10^{-6}$ & $-7.28044 \!\!\times\!\! 10^{-6}$ & $-1.41952 \!\!\times\!\! 10^{-5}$ & $-2.04786 \!\!\times\!\! 10^{-5}$\\
			13 & $-5.70567491243401 \!\!\times\!\! 10^{-6}$ & $-6.20237 \!\!\times\!\! 10^{-6}$ & $-1.21543 \!\!\times\!\! 10^{-5}$ & $-1.47102 \!\!\times\!\! 10^{-5}$\\
			14 & $-3.92760748402515 \!\!\times\!\! 10^{-6}$ & $-7.20853 \!\!\times\!\! 10^{-6}$ & $-9.15766 \!\!\times\!\! 10^{-6}$ & $-1.40720 \!\!\times\!\! 10^{-5}$\\
			15 & $-3.07587337200921 \!\!\times\!\! 10^{-6}$ & $-6.60910 \!\!\times\!\! 10^{-6}$ & $-7.89426 \!\!\times\!\! 10^{-6}$ & $-1.16243 \!\!\times\!\! 10^{-5}$\\
			16 & $-2.22911072118656 \!\!\times\!\! 10^{-6}$ & $-2.71744 \!\!\times\!\! 10^{-6}$ & $-6.39118 \!\!\times\!\! 10^{-6}$ & $-9.15196 \!\!\times\!\! 10^{-6}$\\
			17 & $-1.80153685342167 \!\!\times\!\! 10^{-6}$ & $-3.53017 \!\!\times\!\! 10^{-6}$ & $-5.83645 \!\!\times\!\! 10^{-6}$ & $-8.35031 \!\!\times\!\! 10^{-6}$\\
			18 & $-1.35746592383218 \!\!\times\!\! 10^{-6}$ & $-2.75568 \!\!\times\!\! 10^{-6}$ & $-5.29767 \!\!\times\!\! 10^{-6}$ & $-7.09636 \!\!\times\!\! 10^{-6}$\\
			19 & $-1.1240564465\!\!\times\!\! 10^{-6}     $ & $-2.08440 \!\!\times\!\! 10^{-6}$ & $-4.14648 \!\!\times\!\! 10^{-6}$ & $-5.97031 \!\!\times\!\! 10^{-6}$\\
			20 & $-8.7324870029\!\!\times\!\! 10^{-7}     $ & $-1.84441 \!\!\times\!\! 10^{-6}$ & $-3.79117 \!\!\times\!\! 10^{-6}$ & $-5.15228 \!\!\times\!\! 10^{-6}$\\
			21 & $-7.3719611663\!\!\times\!\! 10^{-7}     $ & $-2.22978 \!\!\times\!\! 10^{-6}$ & $-3.50677 \!\!\times\!\! 10^{-6}$ & $-4.37308 \!\!\times\!\! 10^{-6}$\\
			22 &                                     & $-1.37101 \!\!\times\!\! 10^{-6}$ & $-3.15255 \!\!\times\!\! 10^{-6}$ & $-3.79811 \!\!\times\!\! 10^{-6}$\\
			23 &                                     & $-1.13905 \!\!\times\!\! 10^{-6}$ & $-2.44686 \!\!\times\!\! 10^{-6}$ & $-2.59673 \!\!\times\!\! 10^{-6}$\\
			24 &                                     & $-7.88768 \!\!\times\!\! 10^{-7}$ & $-2.33579 \!\!\times\!\! 10^{-6}$ & $-1.30902 \!\!\times\!\! 10^{-6}$\\
			25 &                                     & $-7.94679 \!\!\times\!\! 10^{-7}$ & $-2.05273 \!\!\times\!\! 10^{-6}$ & $-1.37928 \!\!\times\!\! 10^{-6}$\\
			26 &                                     & $-1.25172 \!\!\times\!\! 10^{-6}$ & $-1.77769 \!\!\times\!\! 10^{-6}$ & $-1.02907 \!\!\times\!\! 10^{-6}$\\
			27 &                                     & $-8.21065 \!\!\times\!\! 10^{-7}$ & $-1.64714 \!\!\times\!\! 10^{-6}$ & $-9.20903 \!\!\times\!\! 10^{-7}$\\
			28 &                                     & $-5.43450 \!\!\times\!\! 10^{-7}$ & $-1.54150 \!\!\times\!\! 10^{-6}$ & $-2.92919 \!\!\times\!\! 10^{-7}$\\
			29 &                                     & $-5.50871 \!\!\times\!\! 10^{-7}$ & $-1.57026 \!\!\times\!\! 10^{-6}$ & $-2.53862 \!\!\times\!\! 10^{-7}$\\
			30 &                                     & $+1.84580 \!\!\times\!\! 10^{-7}$ & $-7.74594 \!\!\times\!\! 10^{-7}$ & $-1.08052 \!\!\times\!\! 10^{-7}$\\
			\hline                                                                   
		\end{tabular}
	\end{table}
}

\begin{figure}[h]
	\centering
	{\includegraphics[width=0.6\linewidth]{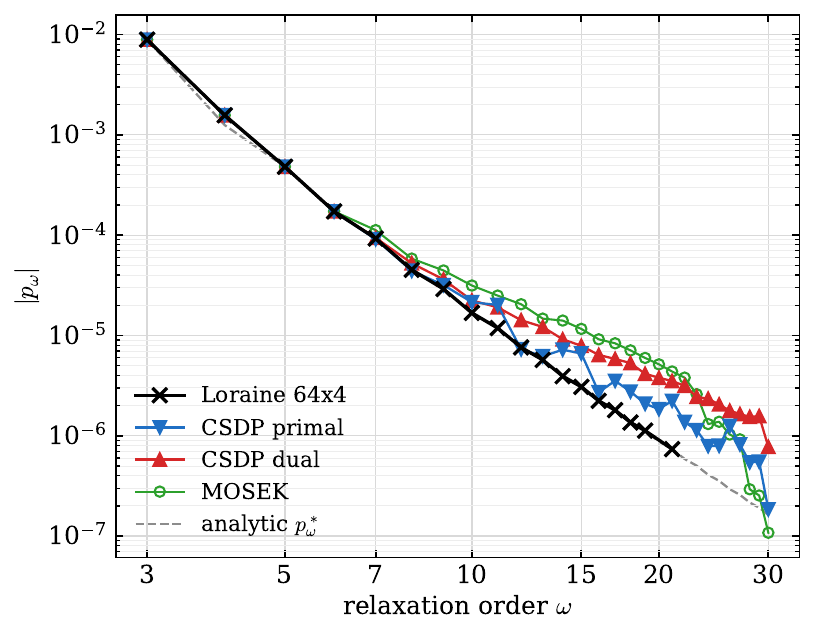}}
	\caption{Problem~2: computed objective $|p_\omega|$ versus relaxation order
		(log-log).  Series as in the legend; the dashed grey line is the reference
		value $p^*_\omega$, extended past the last tabulated \mbox{Loraine\,64x4}
		point to show the reference decay.}
	\label{fig:prob2}
\end{figure}

Starting from $\omega\approx7$, the CSDP and MOSEK values diverge from those computed by
\mbox{Loraine\,64x4}, which follow formula~\cref{eq:p3formula} asymptotically, with exception of orders 4 and 6; see~\Cref{fig:prob2}. 
The CSDP dual sequence remains smooth and almost
monotone; a fit over $\omega\in[7,30]$ has slope $-2.85$ ($R^2=0.989$). On the log-log plot in \Cref{fig:prob2}, this is a
convincingly straight line; yet the conjectured relaxation values~\cref{eq:p3formula}
(confirmed numerically by \mbox{Loraine\,64x4}) decay as $O(\omega^{-4})$, so the CSDP rate is too shallow. MOSEK exhibits a similar shallow-slope false convergence but lies even farther above the high-precision reference curve. This, again, is a sign of false convergence.  The CSDP
primal sequence is increasingly noisy, and both double-precision solvers
eventually lose convergence.

\subsection{Problem 3, stagnation}\label{sec:ex5a}

Finally, we consider the polynomial optimization problem in~$\mathbb{R}^2$:
\begin{equation}\label{eq:p2}
\min_{x\in\mathbb{R}^2}\ x_1
\quad \text{subject to} \quad
x_1^3 - x_2^2\ge 0,\quad
1-x_1^2-x_2^2\geq 0\,.
\end{equation}
As in Problem~1, LICQ is not satisfied at the optimal point $x^*=(0,0)$, due to the construction of the first constraint.
The optimal values of the relaxations were conjectured in
\cite{baldi-kocvara} to be
\begin{equation}\label{eq:p2formula}
f_\omega = -\frac{1}{(2\omega-3)(2\omega-1)}
\end{equation}
and we use this formula as a reference.

\Cref{tab:prob3obj} compares the objective values computed by
\mbox{Loraine\,64x4} with those computed by CSDP and MOSEK. 
As in the previous two problems, the results of \mbox{Loraine\,64x4} agree with the reference rate of \cref{eq:p2formula} to all digits reported in the table.
Starting from
relaxation order~7, the values obtained by the double-precision solvers
deviate increasingly from the reference ones.
{\small
	\begin{table}[h]
		\centering
				\caption{Problem~3, objective values computed by \mbox{Loraine\,64x4}, CSDP and MOSEK}
		\label{tab:prob3obj}
		\begin{tabular}{rlccc}
			\hline
			$\omega$&
			\mbox{Loraine\,64x4}&	CSDP primal & CSDP dual &	MOSEK\\
			\hline
			2 & $-3.33333333333331 \!\!\times\!\! 10^{-1}$  & $-3.33333 \!\!\times\!\! 10^{-1}$ & $-3.33333 \!\!\times\!\! 10^{-1}$ & $-3.33333 \!\!\times\!\! 10^{-1}$\\
			3 & $-6.66666666666658 \!\!\times\!\! 10^{-2}$  & $-6.66666 \!\!\times\!\! 10^{-2}$ & $-6.66667 \!\!\times\!\! 10^{-2}$ & $-6.66666 \!\!\times\!\! 10^{-2}$\\
			4 & $-2.85714285714281 \!\!\times\!\! 10^{-2}$  & $-2.85714 \!\!\times\!\! 10^{-2}$ & $-2.85723 \!\!\times\!\! 10^{-2}$ & $-2.85715 \!\!\times\!\! 10^{-2}$\\
			5 & $-1.58730158730159 \!\!\times\!\! 10^{-2}$  & $-1.58802 \!\!\times\!\! 10^{-2}$ & $-1.70707 \!\!\times\!\! 10^{-2}$ & $-1.58730 \!\!\times\!\! 10^{-2}$\\
			6 & $-1.01010101010101 \!\!\times\!\! 10^{-2}$  & $-1.52905 \!\!\times\!\! 10^{-2}$ & $-1.71158 \!\!\times\!\! 10^{-2}$ & $-1.31161 \!\!\times\!\! 10^{-2}$\\
			7 & $-6.99300699300699 \!\!\times\!\! 10^{-3}$  & $-1.38416 \!\!\times\!\! 10^{-2}$ & $-1.58642 \!\!\times\!\! 10^{-2}$ & $-1.13921 \!\!\times\!\! 10^{-2}$\\
			8 & $-5.12820512820512 \!\!\times\!\! 10^{-3}$  & $-1.60352 \!\!\times\!\! 10^{-2}$ & $-1.84893 \!\!\times\!\! 10^{-2}$ & $-1.19818 \!\!\times\!\! 10^{-2}$\\
			9 & $-3.92156862745097 \!\!\times\!\! 10^{-3}$  & $-1.39454 \!\!\times\!\! 10^{-2}$ & $-1.60108 \!\!\times\!\! 10^{-2}$ & $-1.22985 \!\!\times\!\! 10^{-2}$\\
			10 & $-3.09597523219812 \!\!\times\!\! 10^{-3}$  & $-1.37857 \!\!\times\!\! 10^{-2}$ & $-1.56739 \!\!\times\!\! 10^{-2}$ & $-1.15962 \!\!\times\!\! 10^{-2}$\\
			11 & $-2.50626566416038 \!\!\times\!\! 10^{-3}$  & $-1.54512 \!\!\times\!\! 10^{-2}$ & $-1.78450 \!\!\times\!\! 10^{-2}$ & $-1.24587 \!\!\times\!\! 10^{-2}$\\
			12 & $-2.07039337474118 \!\!\times\!\! 10^{-3}$  & $-1.49811 \!\!\times\!\! 10^{-2}$ & $-1.71769 \!\!\times\!\! 10^{-2}$ & $-1.17524 \!\!\times\!\! 10^{-2}$\\
			13 & $-1.73913043478261 \!\!\times\!\! 10^{-3}$  & $-1.39618 \!\!\times\!\! 10^{-2}$ & $-1.59746 \!\!\times\!\! 10^{-2}$ & $-1.26436 \!\!\times\!\! 10^{-2}$\\
			14 & $-1.48148148148148 \!\!\times\!\! 10^{-3}$  & $-1.43422 \!\!\times\!\! 10^{-2}$ & $-1.64529 \!\!\times\!\! 10^{-2}$ & $-1.10675 \!\!\times\!\! 10^{-2}$\\
			15 & $-1.27713920817166 \!\!\times\!\! 10^{-3}$  & $-1.48284 \!\!\times\!\! 10^{-2}$ & $-1.70557 \!\!\times\!\! 10^{-2}$ & $-1.23824 \!\!\times\!\! 10^{-2}$\\
			\hline
		\end{tabular} 
	\end{table}
}%

\begin{figure}[h]
	\centering
	{\includegraphics[width=0.6\linewidth]{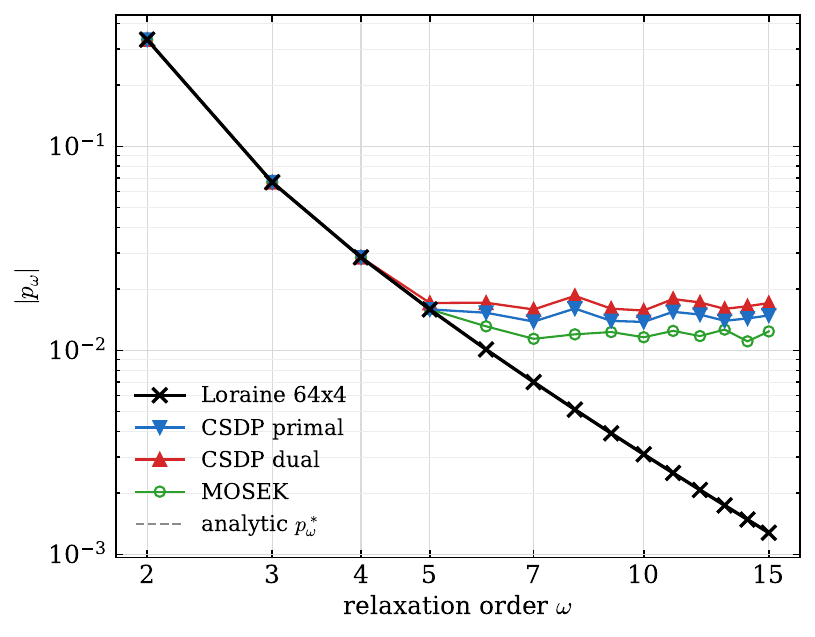}}
	\caption{Problem~3: computed objective $|f_\omega|$ versus relaxation order
		(log-log).  Series as in the legend.}
	\label{fig:prob3stag}
\end{figure}

The resulting objective values as a function of the relaxation order are shown
in \Cref{fig:prob3stag} on a log-log scale. While
\mbox{Loraine\,64x4} follows formula~\cref{eq:p2formula}, the other
solvers stagnate on plateaus of order $10^{-2}$ from
$\omega=7$ onwards.

\section{Spectral and numerical diagnostics}\label{sec:analysis}

\subsection{Positive spectral parts}\label{sec:positive-parts}

Let $X$ and $S$ denote the computed principal paired blocks from a Loraine
reference run; i.e., the moment matrix ($S$) and its dual Gram matrix. For a positive semidefinite matrix $M$ with eigendecomposition
$M=Q\Lambda Q^T$, we define its numerically positive part by retaining the
eigenvalues above the relative cutoff
\begin{equation}\label{eq:cutoff}
 \tau(M)=10^{-12}\lambda_{\max}(M).
\end{equation}
Let $Q_+$ contains the corresponding eigenvectors; then we define
\begin{equation}\label{eq:positive-part}
 M_+=Q_+^TMQ_+,
 \qquad
 \kappa(M_+)=\frac{\lambda_{\max}(M_+)}{\lambda_{\min}(M_+)}.
\end{equation}
The notation $M_+$ denotes a restriction to a numerically selected
subspace, not the usual projection of an indefinite matrix onto the PSD
cone.

Obviously, the full condition numbers of $X$ and $S$ are not useful for the
purpose of measuring conditioning of the SDP problem. Along the central path, complementarity forces certain eigenvalues of $X$ and $S$
toward zero and so the full
condition numbers of both matrices diverge as the barrier parameter decreases, even if the
nonzero spectral clusters of the exact solution are well resolved. We
therefore examine condition numbers of the restricted matrices $X_+$ and~$S_+$.

This requires certain care as the numerical cutoff \cref{eq:cutoff} does not always determine an exact rank. We use
the spectra obtained by \mbox{Loraine\,64x4} with DIMACS tolerance $10^{-16}$ because the positive and near-zero clusters are
visibly separated there and agree with still tighter
runs in selected test cases. Changing the cutoff tolerance in \cref{eq:cutoff} moderately leaves the reported ranks
unchanged. 

For a fixed relaxation order, if the positive cluster of the exact solution is separated
from zero and the computed matrix is sufficiently close in operator norm,
standard invariant-subspace perturbation results imply stability of that
cluster relative to the spectral gap~\cite[Thm.~8.1.10]{GolubVanLoan}.
This supports the use of positive parts once two clusters are clearly
resolved. It does not give a uniform guarantee along the hierarchy, because
the gap separating the positive eigenvalues of the exact solution from zero may shrink with $\omega$.

\subsection{Numerical ranks and stopping tolerance}\label{sec:tolerance}

The complete eigenvalue data for our three problems are given in
Appendix~\ref{app:spectra}; the main text uses only the quantities
needed for the analysis. The computed spectra show that the routinely
used DIMACS error tolerance need not ensure a reliable numerical rank.
For Problem~2 at $\omega=9$, for example, the numerical ranks of the
computed $X$ and $S$ blocks change from $51$ and $48$ at tolerance
$10^{-6}$ to $12$ and $43$ at tolerance $10^{-16}$. For Problem~3 at
$\omega=6$, they change from $26$ and $18$ to $22$ and $6$.

In Problem 2, the numerical rank of \(S\) is more stable under tolerance tightening than that of \(X\). In Problem 3, the rank of \(S\) changes more substantially, while the rank of \(X\) changes by four at each of the three orders considered. Several eigenvalues
classified as positive in the lower-accuracy primal block move into
the near-zero cluster in the tighter computation. Thus, numerical
ranks inferred from a moderately accurate iterate can be misleading
even when feasibility, gap, and complementarity residua appear
satisfactory.

At the tighter tolerance, the numerical ranks satisfy
\begin{equation}\label{eq:numerical-sc}
\operatorname{rank}_{\mathrm{num}}(X)
+\operatorname{rank}_{\mathrm{num}}(S)=m,
\end{equation}
for the selected paired block; see \Cref{tab:ranks}. This is
numerically consistent with strict complementarity (although it does
not prove it).

\begin{table}[h]
	\centering
	\caption{Numerical ranks of the computed $X$ and $S$ blocks
		for $\varepsilon_{_{\mathrm{DIMACS}}}=10^{-16}$.}
	\label{tab:ranks}
	\begin{tabular}{ccccc}
		\toprule
        Problem & $\omega$ & $m$ &
$\operatorname{rank}_{\mathrm{num}}(X)$ &
$\operatorname{rank}_{\mathrm{num}}(S)$ \\
		\midrule
		1 & 5  & 6  & 1 & 5  \\
		1 & 6  & 7  & 1 & 6  \\
		1 & 7  & 8  & 1 & 7  \\
		1 & 8  & 9  & 1 & 8  \\
		1 & 9  & 10& 1 & 9 \\
		\midrule
		2 & 6 & 28 & 1  & 27 \\
		2 & 7 & 36 & 5  & 31 \\
		2 & 8 & 45 & 7  & 38 \\
		2 & 9 & 55 & 12& 43 \\
		\midrule
		3 & 4 & 15 & 11  & 4 \\
		3 & 5 & 21 & 16  & 5 \\
		3 & 6 & 28 & 22  & 6 \\
		\bottomrule
	\end{tabular}
\end{table}

\subsection{Numerical sensitivity of the Nesterov-Todd scaling}\label{sec:nt}

In a primal-dual interior-point method, the Newton system couples the primal and dual semidefinite matrices through the linearized complementarity conditions. The Nesterov-Todd (NT) direction provides a symmetric primal-dual scaling for these conditions.

Although we focus on NT scaling, the numerical difficulties associated with constructing the scaling and solving the resulting Newton system are not specific to this choice; other symmetric primal-dual scalings, such as the Helmberg-Kojima-Monteiro (HKM) scaling, also depend on the current primal and dual matrices and can exhibit similar sensitivity.

The NT scaling matrix $W\succ0$ is characterized by
\begin{equation}\label{eq:nt-balance}
 WSW=X.
\end{equation}
Equivalently,
\begin{equation}\label{eq:nt-W}
 W=X^{1/2}(X^{1/2}SX^{1/2})^{-1/2}X^{1/2}.
\end{equation}
This scaling maps the primal and dual matrices to the same balanced matrix.
Each Newton step of the IPM requires solving the
{Schur complement system}
\begin{equation}\label{eq:schur}
 M\,\Delta y=h,
 \qquad
 M_{ij}=\langle A_i,WA_jW\rangle.
\end{equation}
The important point is that
$W$, and hence the Newton matrix $M$, depends on both $X$ and $S$. 

In Loraine, the NT scaling matrix $W$ is computed using Cholesky factorizations of the primal and dual matrices, followed by an SVD:
\begin{equation*}
X=L_XL_X^T,
\qquad S=L_SL_S^T,
\end{equation*}
\begin{equation*}
L_S^TL_X=U\Sigma V^T.
\end{equation*}
The scaling matrix is then obtained as
\begin{equation*}
G=L_XV\Sigma^{-1/2},
\qquad W=GG^T.
\end{equation*}
Thus, the computation of $W$ depends on both Cholesky factorizations and on the singular values in $\Sigma$. In particular, small singular values are amplified by the inverse square root, making the accuracy of the computed scaling sensitive to the conditioning of the primal and dual matrices.
Other symmetric directions, such as HKM, lead to different formulas but the
same broad conclusion: the Newton equations inherit sensitivity from the
current primal-dual matrices as well as from the problem data.

It is useful to distinguish two computational stages:
\begin{description}
	\item[Stage 1:] the
	scaling matrix $W$ is constructed and the Schur matrix $M$ is assembled;
	\item[Stage 2:] the Newton system \cref{eq:schur} is solved.
\end{description} 
Numerical error can enter at both
stages. Both depend on the primal and dual iterates, while the conditioning
of the second stage also depends essentially on the data matrices $A_i$, their
scaling, and their orientation relative to the eigenspaces of $W$. 

Numerical safeguards such as regularization of a nearly singular Schur
matrix and iterative refinement of the Newton direction address errors
in Stage~2. They cannot, however, correct errors already introduced in
Stage~1. In particular, if the conditioning of the positive spectral
parts of $X$ and $S$ compromises the computation of $W$, the resulting
Schur matrix $M$ may be inaccurate, and solving the corresponding
linear system accurately may still yield an inaccurate Newton
direction.

\subsection{Effective condition number}\label{sec:keff}

The preceding discussion motivates the definition of the \emph{effective condition number}
\begin{equation}\label{eq:keff}
 \kappa_{\mathrm{eff}}:=\kappa(X_+)\kappa(S_+).
\end{equation}
We use $\kappa_{\mathrm{eff}}$ as an empirical indicator of the
combined conditioning of the selected primal and dual blocks. It does
not account for the conditioning of the Schur system arising from the
constraint matrices $A_i$ and their interaction with the scaling matrix
$W$. It also omits the effect of the relative orientation of the relevant
spectral subspaces, the barrier parameter, the remaining PSD blocks,
and implementation-specific error propagation.

The effective condition number can also be expressed in terms of the approximate number of decimal digits potentially lost to conditioning. For arithmetic carrying approximately $p$ decimal digits, this motivates the heuristic precision budget
\begin{equation*}
d_{\mathrm{heur}}=p-\log_{10}\kappa(X_+)-\log_{10}\kappa(S_+)
=p-\log_{10}\kappa_{\mathrm{eff}}.
\end{equation*}
This expression estimates the number of decimal digits remaining after
accounting for the conditioning represented by $\kappa_{\mathrm{eff}}$.
It is not an estimate of the actual number of correct digits in the
computed solution or objective, but a convenient way to compare the
precision demands of different relaxation orders.

For example, in double precision, with approximately $p=16$ decimal
digits, a target of seven digits corresponds to
$\kappa_{\mathrm{eff}}\approx 10^9$. The values in
\Cref{tab:diagnostics} illustrate how rapidly the effective condition
number can grow with the relaxation order. In Problems~1 and~2,
substantial increases in $\kappa_{\mathrm{eff}}$ accompany increases
in the relative error of the computed objective. Problem~3, however,
shows that the relationship is not a sharp threshold: at order~4,
the relative error is small despite $\kappa_{\mathrm{eff}}$ exceeding
$10^8$, whereas at order~5 the error increases substantially.
Thus, $\kappa_{\mathrm{eff}}$ is useful for assessing the precision
demand of a relaxation, but it does not by itself predict the accuracy
of the computed objective. The fourth
column of \Cref{tab:diagnostics} reports the spectral imbalance ratio
$r$, discussed separately in \Cref{sec:imbalance}, which helps
distinguish the strongly dual-dominated cases from the case in which
the primal and dual condition numbers are comparable.

\begin{table}[ht]
\centering
\caption{Effective condition number, spectral imbalance, and
relative error of the CSDP dual objective.}\label{tab:diagnostics}
\begin{tabular}{ccccc}
\toprule
Problem & $\omega$ & $\kappa_{\mathrm{eff}}$ &
$r=\kappa(S_+)/\kappa(X_+)$ & relative error\\
\midrule
1&5&$1.09\times10^4$&$1.09\times10^4$&$2.0\times10^{-6}$\\
1&6&$1.52\times10^5$&$1.52\times10^5$&$1.5\times10^{-5}$\\
1&7&$1.51\times10^6$&$1.51\times10^6$&$1.8\times10^{-3}$\\
1&8&$1.65\times10^7$&$1.65\times10^7$&$1.8\times10^{-2}$\\
1&9&$1.40\times10^8$&$1.40\times10^8$&$1.4\times10^{-1}$\\
\midrule
2&6&$3.72\times10^7$&$3.72\times10^7$&$5.0\times10^{-5}$\\
2&7&$9.19\times10^9$&$1.15\times10^6$&$1.5\times10^{-2}$\\
2&8&$4.34\times10^{13}$&$1.57\times10^6$&$1.6\times10^{-1}$\\
2&9&$2.22\times10^{14}$&$1.08\times10^6$&$2.4\times10^{-1}$\\
\midrule
3&4&$5.39\times10^8$&$2.06$&$3.0\times10^{-5}$\\
3&5&$2.04\times10^{12}$&$1.96$&$7.8\times10^{-2}$\\
3&6&$5.72\times10^{15}$&$1.86$&$7.0\times10^{-1}$\\
\bottomrule
\end{tabular}
\end{table}

Problem~2 further illustrates that the product
$\kappa_{\mathrm{eff}}$ need not track the objective error as closely
as either of its factors. At the largest reported orders, the absolute
error in the CSDP dual objective follows the dual spectral factor more
closely than the full product. This is consistent with the strongly
dual-dominated spectrum in this example. We therefore use
$\kappa_{\mathrm{eff}}$ as a combined warning indicator of precision
demand, rather than as a formula for the objective error.

\subsection{Objective value sensitivity}\label{sec:objective}

For problems with a small objective value, possibly tending to zero as the relaxation order increases, even a modest error in the computed objective can lead to a large relative error.

For a computed objective value $f_\omega$ and a trusted reference
$f_\omega^*\ne0$, we define the normalized objective error
\begin{equation*}
\eta_{\mathrm{obj}}
:=\frac{|f_\omega-f^*_\omega|}
{\|B\|_F\;\lambda_{\max}(X_+)},
\end{equation*}
and the amplification factor
\begin{equation}\label{eq:C}
C:=\frac{\|B\|_F\;\lambda_{\max}(X_+)}{|f^*_\omega|},
\end{equation}
where $X_+$ is the positive spectral part of the computed primal
solution. 
Then
\begin{equation*}
\frac{|f_\omega-f^*_\omega|}{|f^*_\omega|}
=\eta_{\mathrm{obj}}\,C.
\end{equation*}
Thus, $C$ measures how the normalized objective error is amplified
when expressed as a relative error in the objective. 


The scale in the definition of $C$ could instead be based on the full block-diagonal primal matrix and objective, including the localizing matrices and their duals. However, the available spectral tables contain only the principal paired block $(X,S)$; accordingly, \Cref{tab:objcond} uses $\|B\|_F\;\lambda_{\max}(X_+)$. The resulting factor $C$ is therefore a normalization diagnostic for this block and this SDP representation, rather than a condition measure for the full SDP.

\begin{table}[ht]
	\centering
	\caption{Principal-block objective amplification factor $C$, with
		$\|B\|_F=1$. Values of $\lambda_{\max}(X_+)$ are those of the
		$\varepsilon_{\mathrm{DIMACS}}=10^{-16}$ runs.}
	\label{tab:objcond}
	\begin{tabular}{ccccc}
		\toprule
		Problem & $\omega$ & $\lambda_{\max}(X_+)$ & $|f_\omega|$ & $C$\\
		\midrule
		1&5&$1.60\times10^2$&$6.67\times10^{-2}$&$2.4\times10^3$\\
		1&6&$9.22\times10^2$&$4.17\times10^{-2}$&$2.2\times10^4$\\
		1&7&$5.25\times10^3$&$2.86\times10^{-2}$&$1.8\times10^5$\\
		1&8&$2.98\times10^4$&$2.08\times10^{-2}$&$1.4\times10^6$\\
		1&9&$1.69\times10^5$&$1.59\times10^{-2}$&$1.1\times10^7$\\
		\midrule
		2&6&$1.26\times10^2$&$1.72\times10^{-4}$&$7.3\times10^5$\\
		2&7&$4.03\times10^2$&$9.26\times10^{-5}$&$4.4\times10^6$\\
		2&8&$3.30\times10^3$&$4.51\times10^{-5}$&$7.3\times10^7$\\
		2&9&$9.89\times10^3$&$2.91\times10^{-5}$&$3.4\times10^8$\\
		\midrule
		3&4&$2.39\times10^3$&$2.86\times10^{-2}$&$8.4\times10^4$\\
		3&5&$8.17\times10^4$&$1.59\times10^{-2}$&$5.2\times10^6$\\
		3&6&$2.79\times10^6$&$1.01\times10^{-2}$&$2.8\times10^8$\\
		\bottomrule
	\end{tabular}
\end{table}

The table shows that $C$ grows by several orders of magnitude with
the relaxation order in all three examples. This growth reflects the
increasing ratio of the primal scale, measured by
$\lambda_{\max}(X_+)$, to the magnitude of the objective value.
Consequently, an objective error that is small relative to
$\|B\|_F\,\lambda_{\max}(X_+)$ can correspond to a large relative error
in the reported objective. This amplification does not identify the
numerical origin of the objective error. The effective condition number
$\kappa_{\mathrm{eff}}$, in contrast, characterizes the combined
conditioning of the selected primal and dual blocks, but does not
account for amplification by a small objective value. The two
quantities therefore describe different aspects of the observed loss
of accuracy.

\subsection{False convergence, stagnation, and spectral imbalance}
\label{sec:imbalance}

We define the \emph{spectral imbalance ratio}
\begin{equation}\label{eq:ratio}
r:=\frac{\kappa(S_+)}{\kappa(X_+)}.
\end{equation}
A large value of $r$ indicates that the selected dual block is much
more ill-conditioned than the primal block. In our examples, this
imbalance is associated with false convergence, whereas large values of \(\kappa_{\mathrm{eff}}\), with comparable primal and dual condition numbers, are associated with stagnation; see \Cref{tab:diagnostics}.
Thus, $r$ provides a possible indicator of the type of numerical
failure, rather than a criterion that defines or predicts it.

For Problem~1, $X_+$ is one-dimensional and $\kappa(X_+)=1$, while
$\kappa(S_+)$ grows rapidly. Problem~2 is also strongly dual-dominated
over the displayed orders, although its primal block has increasing
rank and its primal objective shows increasing numerical variability. Both problems exhibit
a smooth but incorrect dual sequence. In Problem~3, the two
positive-part condition numbers are comparable and both become large;
the computed hierarchy loses its systematic trend and stagnates.

These observations suggest that a large spectral imbalance ratio may
serve as a warning sign for false convergence, while comparable
primal and dual conditioning may accompany stagnation. However, three
examples are far too few to establish this as a general principle.
In particular, the ratio omits problem data and all other
blocks, and Problem~1 does not separate the effect of a large ratio
from the special fact that $X_+$ has rank one.

\subsection{Cholesky inverse residuum}\label{sec:rho}

The diagnostics discussed above rely on reliable spectral information
about the primal and dual solutions or on a trusted reference objective
value, which may require a high-precision computation. What can be
monitored when only a standard solver is available, without an
independent reference solution? We consider a simple diagnostic that
can be evaluated during the interior-point iteration.

A practical diagnostic available during an interior-point iteration is
\begin{equation}\label{eq:rho}
 \rho:=\|SS^{-1}_{\mathrm{comp}}-I\|,
\end{equation}
where $S^{-1}_{\mathrm{comp}}v$ denotes the inverse action obtained by
solving with the computed Cholesky factors. No explicit inverse need be
formed. The quantity measures the quality of this inverse action and uses
the same factorization involved in constructing the NT scaling.

For backward-stable factorization and triangular solves one expects the
residuum to grow with machine precision and the conditioning of the full
current iterate $S$; see \cite[Thms.~10.3,10.4]{HighamAccuracy}. In the experiments,
loss of useful progress was associated with values around $10^{-7}$ in
Float64 and around $10^{-20}$ in Float64x4. These levels are empirical and
solver-dependent; the latter is also influenced by the $10^{-16}$ stopping
rule.

The residuum is attractive because it is inexpensive and available during
the solve. Its limitation is equally important: it monitors one
factorization-based operation involving $S$ and contains no direct
information about $X$, the constraint operator, or objective amplification.
It is a warning signal, not a universal stopping criterion.

In solvers that support arithmetic at multiple precisions, the residuum $\rho$ could also serve as an early-warning indicator for switching to higher precision. Such a switch should be considered while the Schur matrix is still positive definite, rather than waiting for indefiniteness to occur, by which point numerical difficulties may already have compromised the iteration. The appropriate threshold for $\rho$ would need to be calibrated for the solver and problem class.

\subsection{Practical use of the diagnostics}\label{sec:practical}

\Cref{tab:roles} summarizes the diagnostics in increasing order of availability in IP software, a short checklist for anyone computing moment-SOS bounds in finite precision.%
\newcolumntype{P}[1]{>{\raggedright\arraybackslash}p{#1}}%
\begin{table}[ht]
	\centering
	\caption{Interpretation of the numerical diagnostics.}
	\label{tab:roles}
	\begin{tabular}{P{.17\linewidth}P{.16\linewidth}P{.27\linewidth}P{.27\linewidth}}
		\toprule
		Quantity & Availability & Main use & Main limitation\\
		\midrule
		$\rho$ \cref{eq:rho}
		& during the solve
		& warns of loss of accuracy in a Cholesky-based inverse action
		& omits $X$, constraint geometry, and objective amplification\\
		\midrule
		Primal-dual gap
		& during the solve
		& monitors agreement between the computed primal and dual objectives
		& a small gap does not imply accurate objective values\\
		\midrule
		$C$ \cref{eq:C}
		& with a trusted reference and block scale
		& measures amplification of normalized objective discrepancy
		& does not identify its source\\
		\midrule
		$\kappa_{\mathrm{eff}}$ \cref{eq:keff}
		& after reliable spectral computation
		& summarizes the positive-part spectral spread of the selected pair
		& block- and representation-dependent; not an error bound\\
		\midrule
		$r$ \cref{eq:ratio}
		& after reliable spectral computation
		& describes primal-dual spectral imbalance
		& no established classification role\\
		\bottomrule
	\end{tabular}
\end{table}%

A reasonable workflow is therefore to monitor factorization residua
during the solve, using the primal-dual gap as a measure of agreement
between the computed objectives rather than as a certificate of their
accuracy. One should inspect the hierarchy as a sequence rather than as
isolated SDP runs, and use higher precision or independent references
when the sequence becomes suspicious. A smooth curve is not by itself
reassuring: it may be the signature of false convergence. Spectral
diagnostics should be computed only from iterates accurate enough to
separate positive and near-zero eigenvalue clusters.

\section{Conclusions}\label{sec:conclusion}

The three examples demonstrate that numerical failure of successive SDP
relaxations can produce two distinct patterns. In false
convergence, the computed sequence of objective values remains smooth but follows the wrong
asymptotic rate. In stagnation, it loses a systematic convergence trend and
forms a plateau or oscillates. The first pattern is particularly dangerous
because standard extrapolation can turn a numerical artefact into an
apparently convincing convergence claim.

Two effects should be kept separate. The objective amplification factor $C$ describes how a normalized discrepancy is magnified when the relaxation value approaches zero. Separately, the numerical analysis of the Newton step reveals a sensitivity to the conditioning of both the primal and dual semidefinite matrices. This motivates the effective condition number
$\kappa_{\mathrm{eff}}=\kappa(X_+)\kappa(S_+)$,
which captures the combined spectral conditioning of the primal and dual solution components. Its growth is associated with the observed deterioration in our examples. However, given the limited number of examples studied, it should not be regarded as a general condition number of an SDP.

The tolerance study gives a further warning: a small DIMACS error does not
ensure that numerical ranks and positive-part condition numbers are already
reliable. In Problems~2 and~3, tightening the tolerance changes the numerically identified
primal positive subspace substantially. Spectral diagnostics therefore
require sufficiently accurate iterates to separate positive and near-zero
eigenvalue clusters; satisfying a standard stopping test may not suffice.

The role of the imbalance ratio $r=\kappa(S_+)/\kappa(X_+)$ is to organize the three examples:
the false-convergence cases are strongly dual-dominated, whereas the
stagnation case has comparable primal and dual factors. Although the pattern is consistent across these examples, its generality remains an open question.

The Cholesky inverse residuum $\rho$ provides a cheap in-solver warning
of loss of accuracy in the inverse action used to form the scaling matrix.
In solvers supporting multiple arithmetic precisions, monitoring $\rho$
may help trigger a switch to higher precision before the Schur complement
matrix becomes indefinite. The appropriate threshold is problem- and
solver-dependent, and $\rho$ does not capture every source of numerical
error.

The central practical conclusion is that solver status, a small
primal-dual gap, and a visually smooth sequence do not certify the accuracy
of high-order moment-SOS bounds. The gap measures agreement between the
computed primal and dual objectives, which may both be far from the true
optimum. When the asymptotic behaviour matters, independent formulas,
reliable spectral checks, alternative formulations, or higher-precision
computations may be necessary.

Related approaches include iterative refinement of approximate SDP solutions
\cite{MohammadisiahroudiEtAl2026} and rounding procedures that exploit an
approximate optimal partition to construct an improved primal-dual
solution \cite{MohammadNezhadTerlaky2019}. Such methods provide possible
routes to improving computed solutions, but their effectiveness for the
high-order moment-SOS relaxations studied here remains to be investigated.

\begin{appendices}

\section{Problem 1 with SDP solvers in Julia/JuMP}\label{app:solvers}

Below are final messages of various SDP solvers installed in Julia/JuMP for Problem~1 with $\omega=11$. These span both interior-point solvers (Clarabel, Hypatia, MOSEK, DSDP, CSDP) and first-order operator-splitting solvers (SCS, COSMO); that the latter fail as well indicates that the failure is not confined to the interior-point methods analysed in the main text, and is consistent with the objective conditioning $C$ of \Cref{sec:objective}, a property of the relaxation rather than of any algorithm. DSDP is the exception among the solvers listed, its value agreeing with the exact one to four digits. Recall that, for this problem, the exact optimal value is $-\frac{1}{\omega(\omega-2)}=-\frac{1}{11(11-2)} = -1.\overline{01}\times 10^{-2}$.

\bigskip
\noindent
{\small
	\begin{minipage}[t]{0.499\textwidth}
		\begin{verbatim}
		Clarabel:
		Primal:-1.9070e-02 Dual:-1.9070e-02
		Terminated with status = solved 
		(reduced accuracy)
		Term. status:  ALMOST_OPTIMAL
		Primal status: NEARLY_FEASIBLE_POINT
		Dual status:   NEARLY_FEASIBLE_POINT
		\end{verbatim}
	\end{minipage}\hfill
	\begin{minipage}[t]{0.499\textwidth}
		\begin{verbatim}
		Hypatia:
		Primal:-1.2795e-02 Dual:-1.2832e-02
		optimal solution found; terminating
		Term. status:  OPTIMAL
		Primal status: FEASIBLE_POINT
		Dual status:   FEASIBLE_POINT
		\end{verbatim}
	\end{minipage}
	
	\noindent
	\begin{minipage}[t]{0.499\textwidth}
		\begin{verbatim}
		Mosek:
		Term. status:  SLOW_PROGRESS
		Primal status: UNKNOWN_RESULT_STATUS
		Dual status:   UNKNOWN_RESULT_STATUS
		\end{verbatim}
	\end{minipage}\hfill
	\begin{minipage}[t]{0.499\textwidth}
		\begin{verbatim}
		SCS:
		objective = -0.060854
		Term. status:  OPTIMAL
		Primal status: FEASIBLE_POINT
		Dual status:   FEASIBLE_POINT
		\end{verbatim}
	\end{minipage}
	
	\noindent
	\begin{minipage}[t]{0.499\textwidth}
		\begin{verbatim}
		COSMO:
		Optimal objective: -0.04686
		Term. status:  OPTIMAL
		Primal status: FEASIBLE_POINT
		Dual status:   FEASIBLE_POINT
		\end{verbatim}
	\end{minipage}\hfill
	\begin{minipage}[t]{0.499\textwidth}
		\begin{verbatim}
		DSDP:
		Optimal lower bound: -0.0101080313
		Term. status:  OPTIMAL
		Primal status: FEASIBLE_POINT
		Dual status:   FEASIBLE_POINT
		\end{verbatim}
	\end{minipage}
	
	\noindent
	\begin{minipage}[t]{0.95\textwidth}
		\begin{verbatim}
		CSDP:
		Stuck at edge of primal feasibility, giving up. 
		Partial Success: SDP solved with reduced accuracy
		Termination status: ALMOST_OPTIMAL
		Primal status:      NEARLY_FEASIBLE_POINT
		Dual status:        NEARLY_FEASIBLE_POINT
		\end{verbatim}
	\end{minipage}
}

\section{Eigenvalue data}\label{app:spectra}
\begin{figure}[h]
	\centering
	{\includegraphics[width=0.49\linewidth]{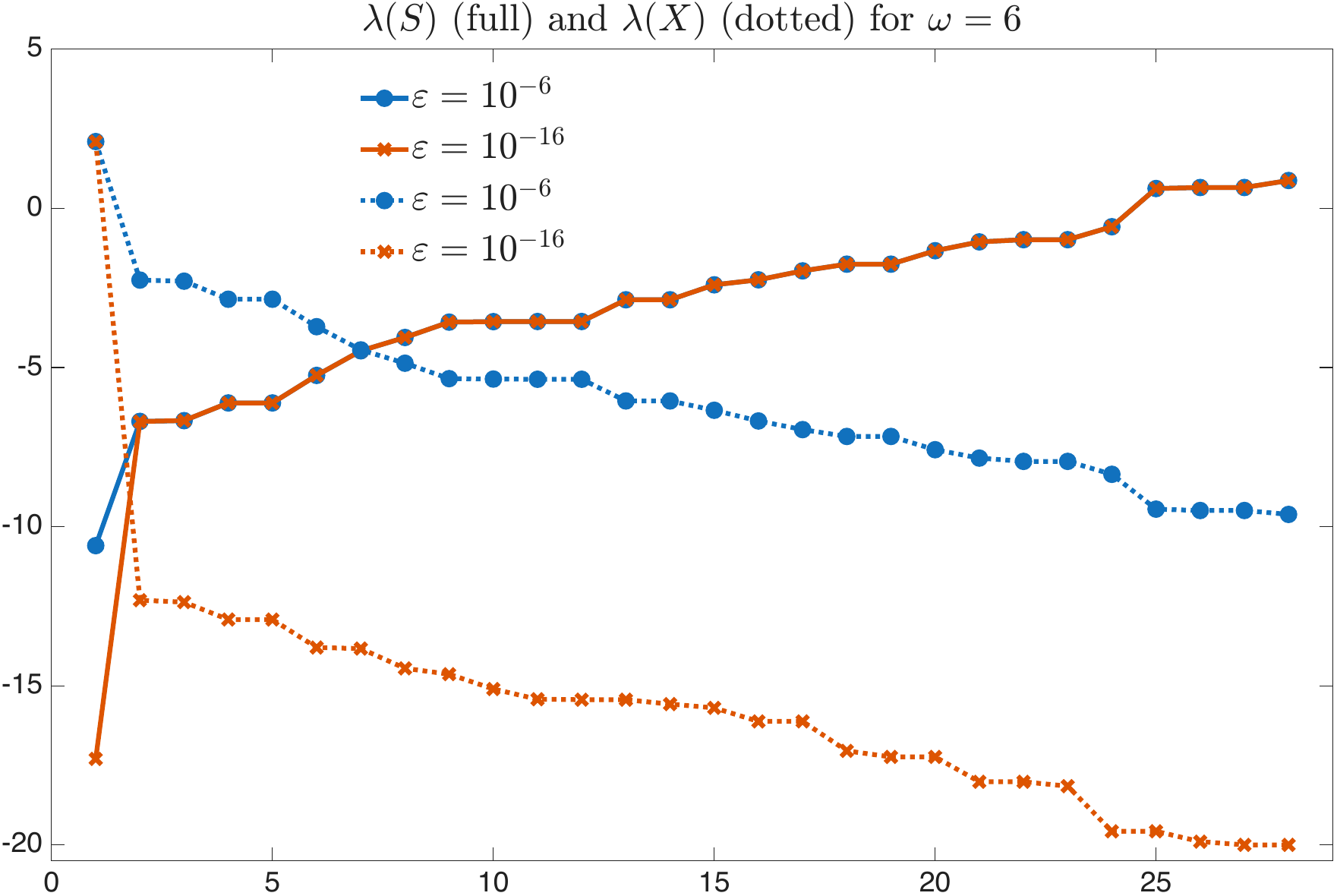}}
	{\includegraphics[width=0.49\linewidth]{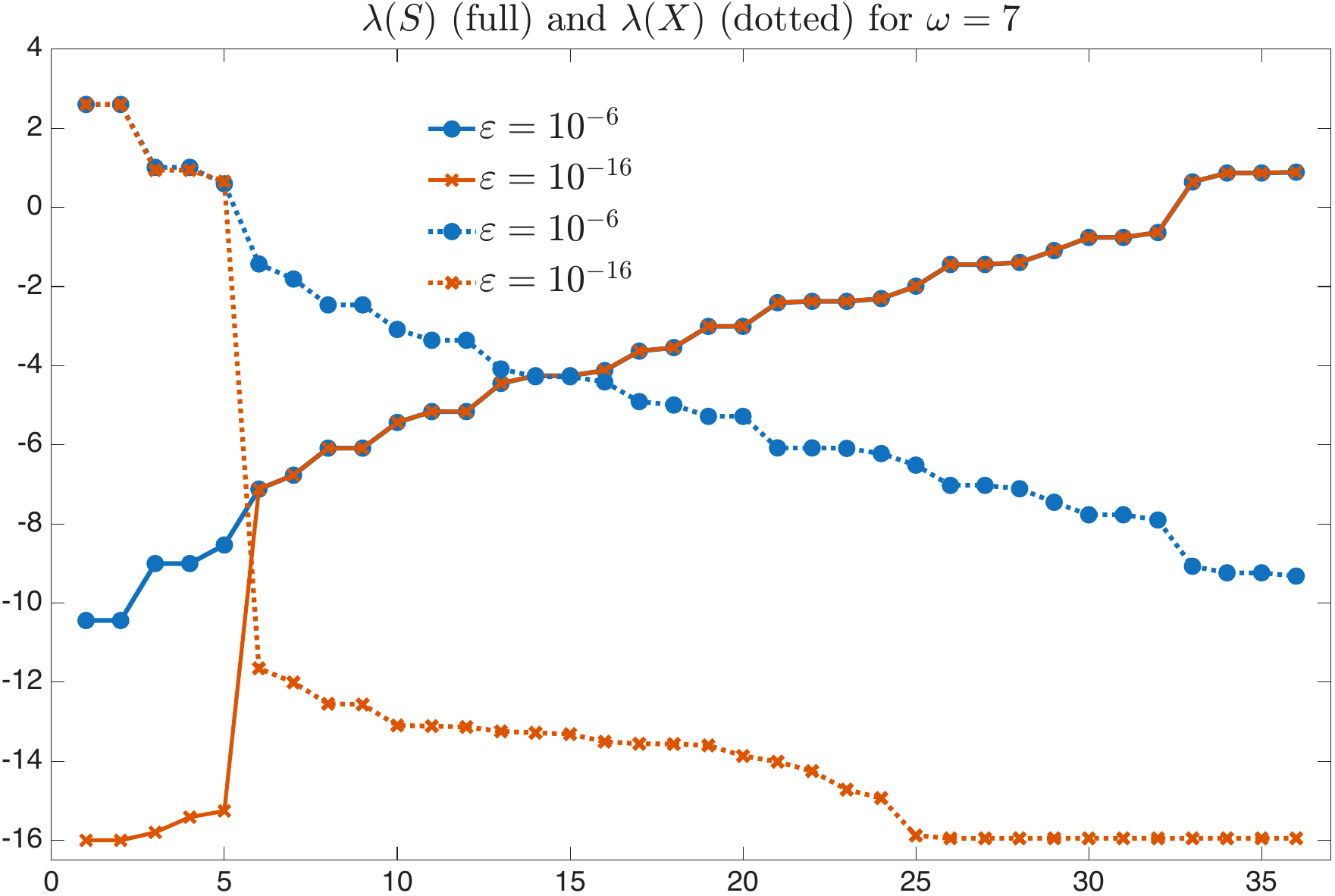}}\\[0.5em]
	{\includegraphics[width=0.49\linewidth]{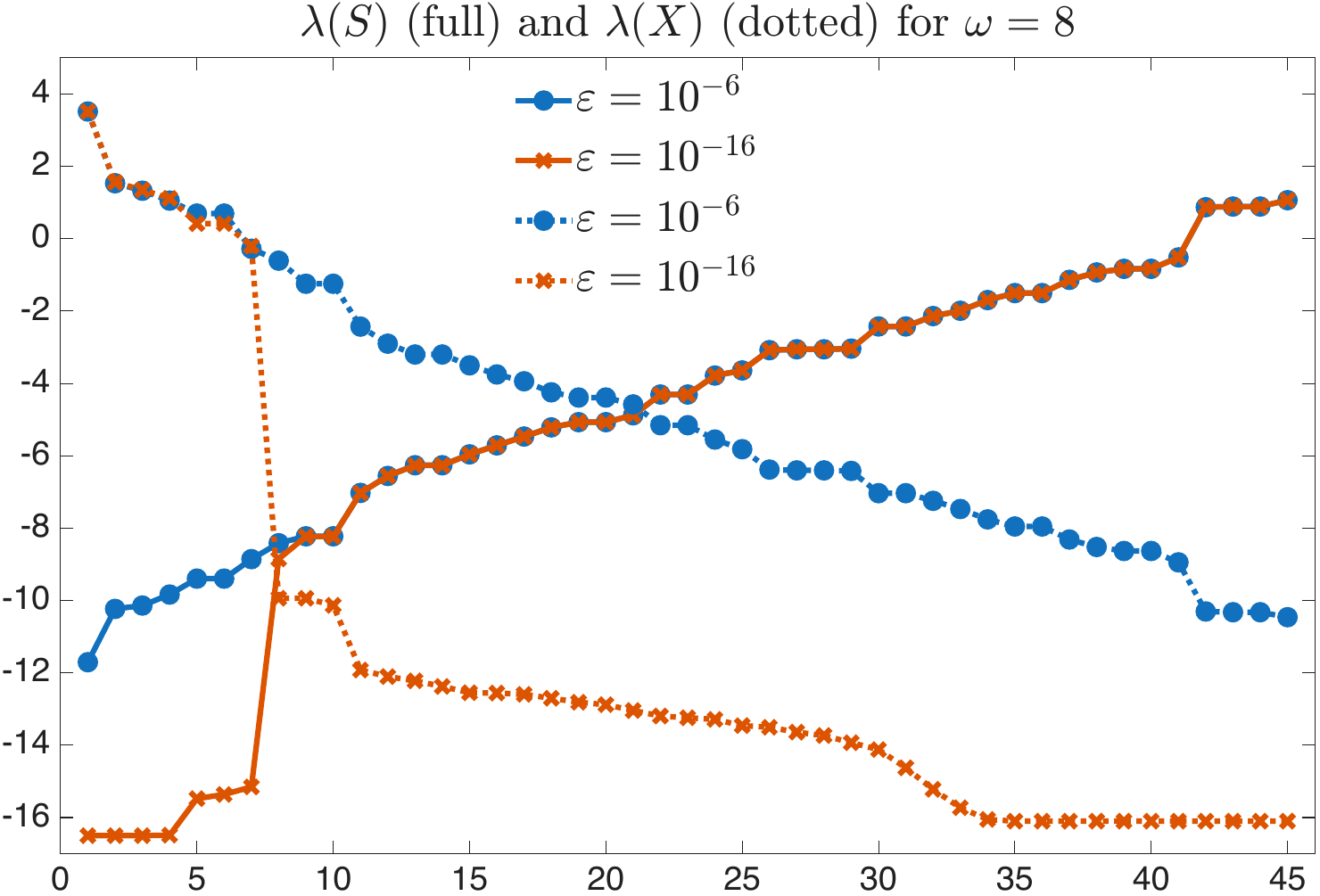}}
	{\includegraphics[width=0.49\linewidth]{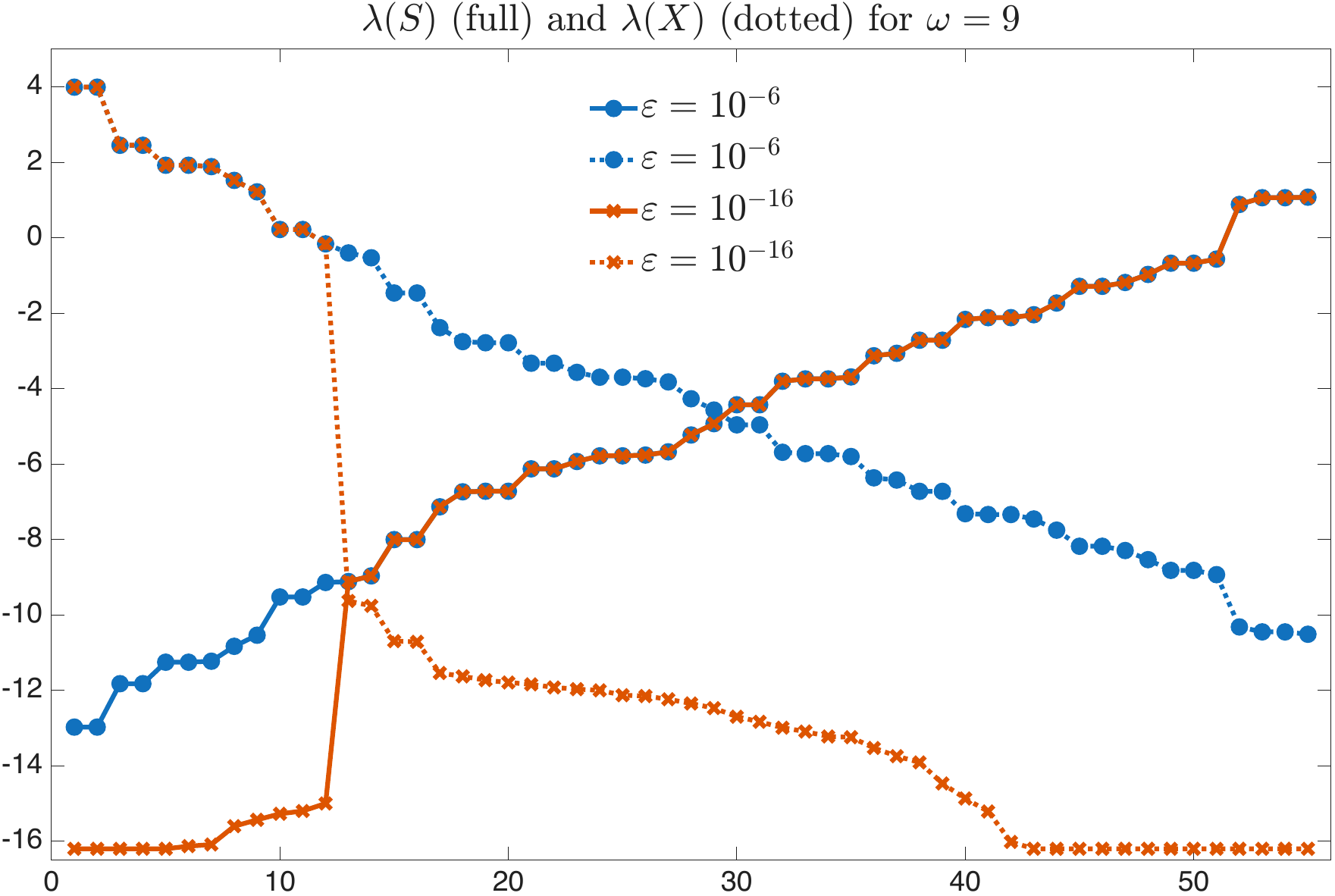}}
	\caption{Problem~2: logarithms of the eigenvalues of $S$ (solid)
		and $X$ (dotted) for $\omega=6,7,8,9$, at
		$\varepsilon_{_{\mathrm{DIMACS}}}=10^{-6}$ (blue) and $10^{-16}$ (red).}
	\label{fig:prob2ev}
\end{figure}

\begin{figure}[h]
	\centering
	{\includegraphics[width=0.49\linewidth]{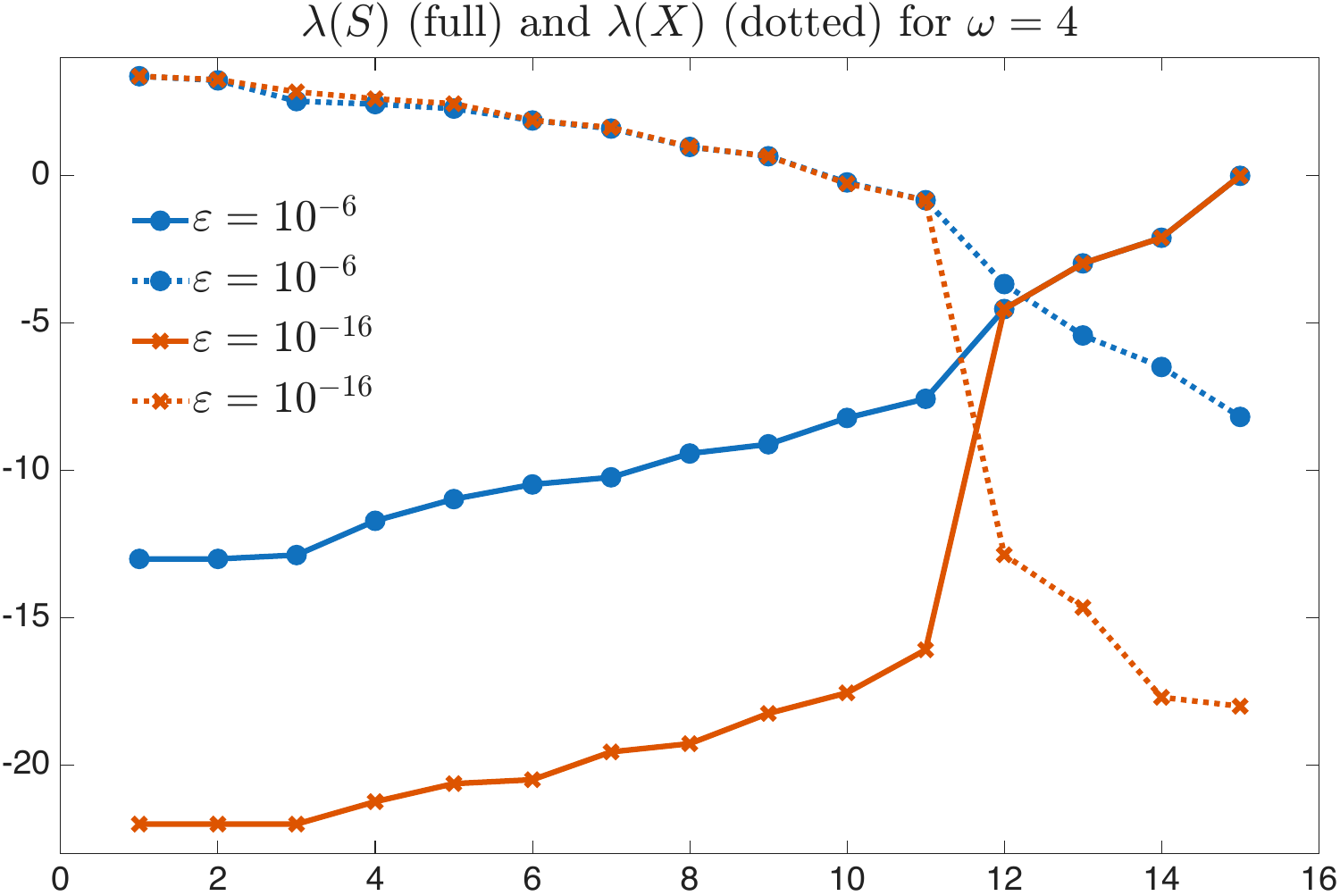}}
	{\includegraphics[width=0.49\linewidth]{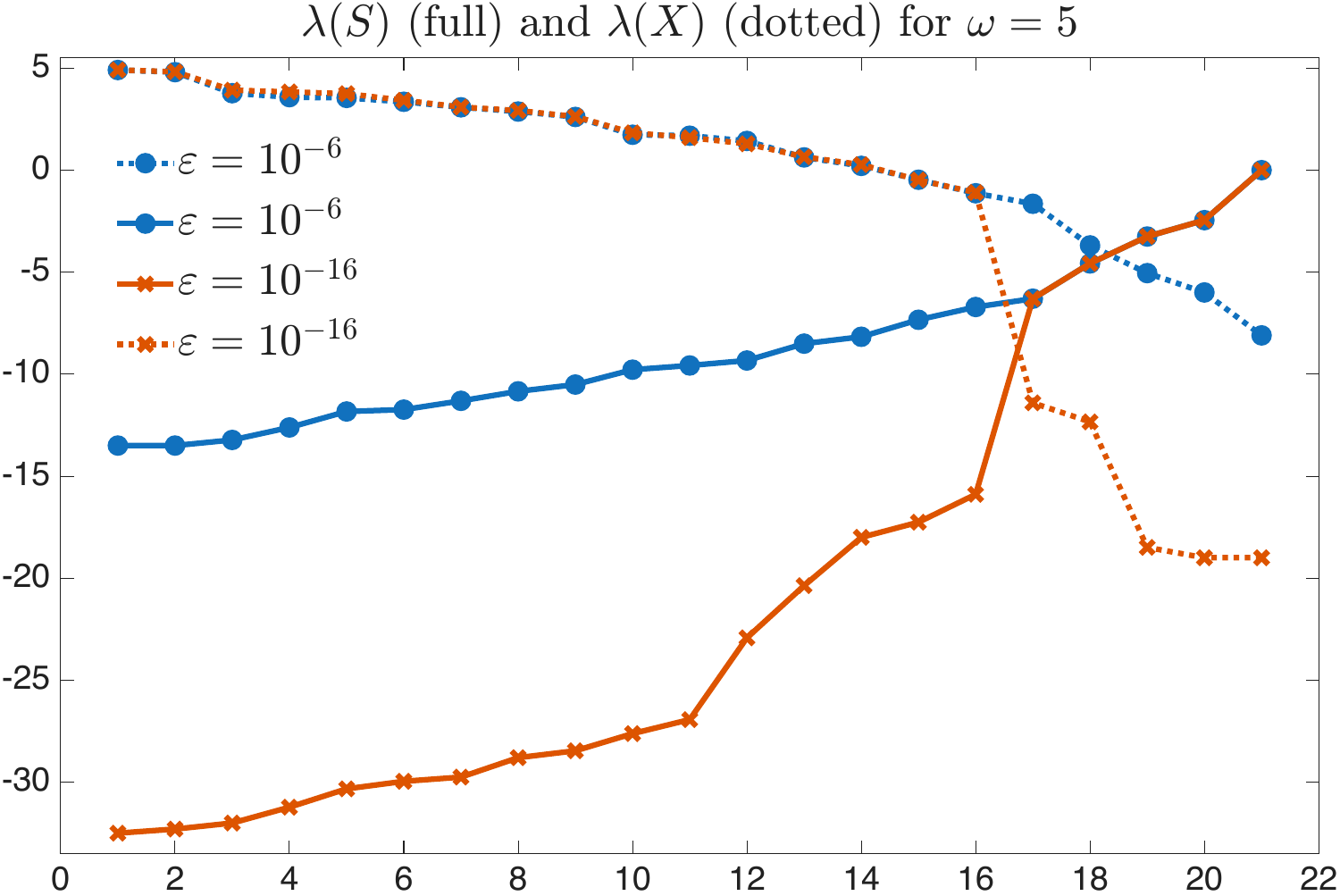}}\\[0.5em]
	{\includegraphics[width=0.49\linewidth]{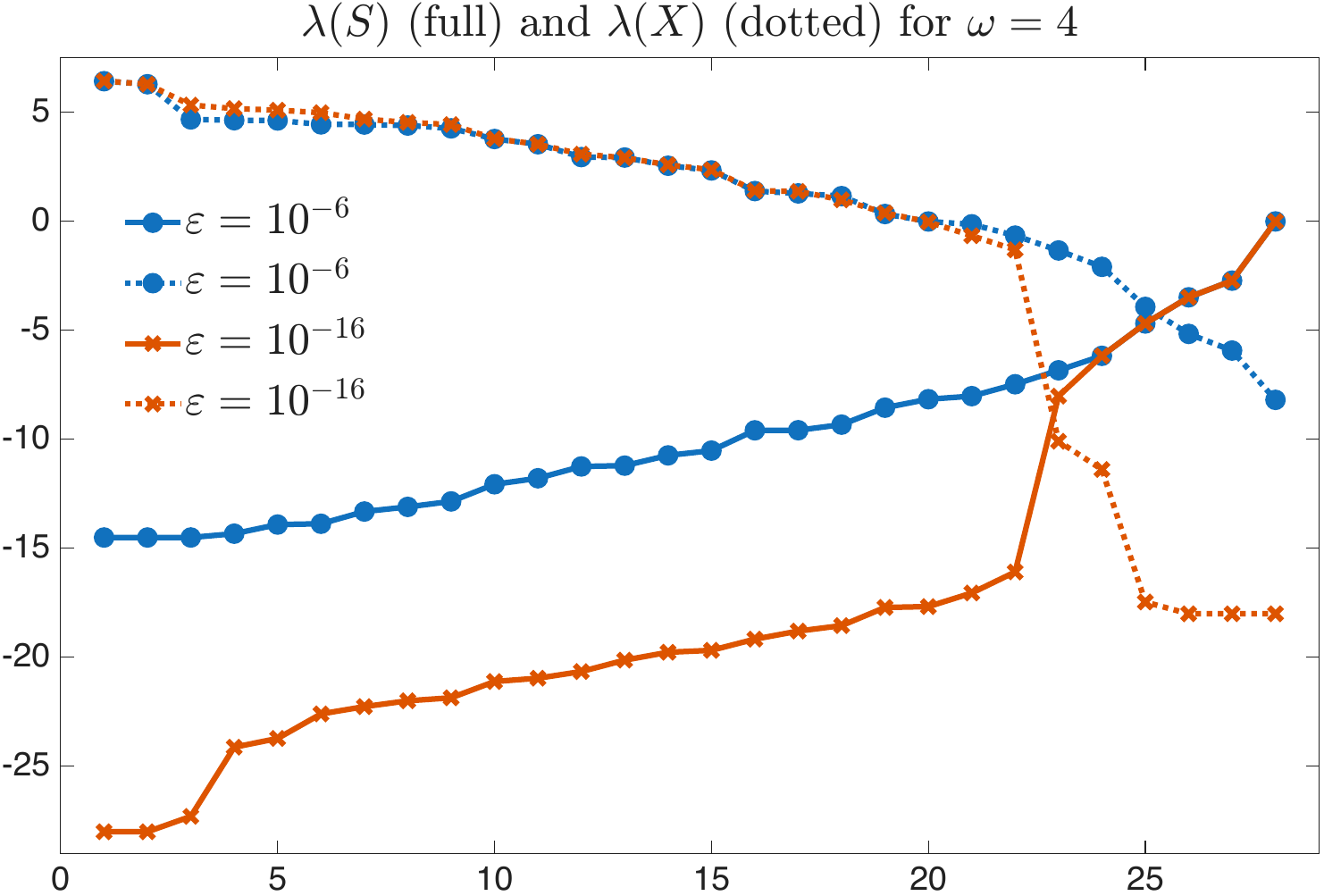}}
	\caption{Problem~3: logarithms of the eigenvalues of $S$ (solid)
		and $X$ (dotted) for $\omega=4,5,6$, at
		$\varepsilon_{_{\mathrm{DIMACS}}}=10^{-6}$ (blue) and $10^{-16}$ (red).}
	\label{fig:prob3ev}
\end{figure}

All eigenvalue data are computed by \mbox{Loraine\,64x4}. For
Problem~1, the eigenvalues of the moment matrix $S$ and the Gram matrix
$X$ are listed explicitly. For Problems~2 and~3, where the matrices are
larger, we report the numerical rank, the range of the positive
eigenvalues, and the corresponding condition number for each stopping
criterion. The values $\lambda_{\max+}$ at
$\varepsilon_{_{\mathrm{DIMACS}}}=10^{-16}$ are those used to form the
objective conditioning $C$ of \Cref{tab:objcond}. \Cref{fig:prob2ev,fig:prob3ev} show the corresponding eigenvalue
distributions on a logarithmic scale.

\subsection*{Problem~1}

\paragraph{Relaxation order $\omega=5$}
\begin{center}
	\begin{tabular}{ll|ll}
		\multicolumn{2}{c|}{$\varepsilon_{_{\mathrm{DIMACS}}}=10^{-6}$} &
		\multicolumn{2}{c}{$\varepsilon_{_{\mathrm{DIMACS}}}=10^{-16}$} \\
		\hline
		$X$ & $S$ & $X$ & $S$ \\
		\hline
		$1.37\!\!\times\!\!10^{-8}$ & $7.47\!\!\times\!\!10^{-10}$ &
		$6.90\!\!\times\!\!10^{-19}$ & $3.62\!\!\times\!\!10^{-20}$ \\
		$2.24\!\!\times\!\!10^{-8}$ & & $9.75\!\!\times\!\!10^{-19}$ & \\
		$1.31\!\!\times\!\!10^{-6}$ & & $6.90\!\!\times\!\!10^{-17}$ & \\
		$4.98\!\!\times\!\!10^{-6}$ & & $2.28\!\!\times\!\!10^{-16}$ & \\
		$6.47\!\!\times\!\!10^{-5}$ & & $3.41\!\!\times\!\!10^{-15}$ & \\
		\hline
		$1.60\!\!\times\!\!10^{2}$  & $3.50\!\!\times\!\!10^{-4}$ &
		$1.60\!\!\times\!\!10^{2}$  & $3.50\!\!\times\!\!10^{-4}$ \\
		& $8.61\!\!\times\!\!10^{-3}$ & & $8.61\!\!\times\!\!10^{-3}$ \\
		& $1.78\!\!\times\!\!10^{-2}$ & & $1.78\!\!\times\!\!10^{-2}$ \\
		& $3.51$              & & $3.51$ \\
		& $3.82$              & & $3.82$ \\
	\end{tabular}
\end{center}

\paragraph{Relaxation order $\omega=6$}
\begin{center}
	\begin{tabular}{ll|ll}
		\multicolumn{2}{c|}{$\varepsilon_{_{\mathrm{DIMACS}}}=10^{-6}$} &
		\multicolumn{2}{c}{$\varepsilon_{_{\mathrm{DIMACS}}}=10^{-16}$} \\
		\hline
		$X$ & $S$ & $X$ & $S$ \\
		\hline
		$1.17\!\!\times\!\!10^{-8}$ & $9.00\!\!\times\!\!10^{-11}$ & $2.93\!\!\times\!\!10^{-13}$ & $2.15\!\!\times\!\!10^{-16}$ \\
		$1.71\!\!\times\!\!10^{-8}$ &  & $1.46\!\!\times\!\!10^{-30}$ &  \\
		$2.36\!\!\times\!\!10^{-6}$ &  & $4.26\!\!\times\!\!10^{-17}$ &  \\
		$6.27\!\!\times\!\!10^{-6}$ &  & $6.00\!\!\times\!\!10^{-17}$ &  \\
		$7.39\!\!\times\!\!10^{-5}$ &  & $1.33\!\!\times\!\!10^{-14}$ &  \\
		$1.25\!\!\times\!\!10^{-3}$ &  & $6.57\!\!\times\!\!10^{-14}$ &  \\
		\hline
		$9.22\!\!\times\!\!10^{2}$ & $3.05\!\!\times\!\!10^{-5}$ & $9.22\!\!\times\!\!10^{2}$ & $3.06\!\!\times\!\!10^{-5}$ \\
		& $5.22\!\!\times\!\!10^{-4}$ &  & $5.22\!\!\times\!\!10^{-4}$ \\
		& $5.18\!\!\times\!\!10^{-3}$ &  & $5.20\!\!\times\!\!10^{-3}$ \\
		& $1.56\!\!\times\!\!10^{-2}$ &  & $1.57\!\!\times\!\!10^{-2}$ \\
		& $3.47\!\!\times\!\!10^{0}$ &  & $3.48\!\!\times\!\!10^{0}$ \\
		& $4.66\!\!\times\!\!10^{0}$ &  & $4.66\!\!\times\!\!10^{0}$ \\
	\end{tabular}
\end{center}

\paragraph{Relaxation order $\omega=7$}
\begin{center}
	\begin{tabular}{ll|ll}
		\multicolumn{2}{c|}{$\varepsilon_{_{\mathrm{DIMACS}}}=10^{-6}$} &
		\multicolumn{2}{c}{$\varepsilon_{_{\mathrm{DIMACS}}}=10^{-16}$} \\
		\hline
		$X$ & $S$ & $X$ & $S$ \\
		\hline
		$1.91\!\!\times\!\!10^{-9}$ & $1.76\!\!\times\!\!10^{-12}$ & $4.80\!\!\times\!\!10^{-13}$ & $7.19\!\!\times\!\!10^{-16}$ \\
		$2.72\!\!\times\!\!10^{-9}$ &  & $1.38\!\!\times\!\!10^{-14}$ &  \\
		$5.48\!\!\times\!\!10^{-7}$ &  & $5.62\!\!\times\!\!10^{-17}$ &  \\
		$1.61\!\!\times\!\!10^{-6}$ &  & $2.20\!\!\times\!\!10^{-15}$ &  \\
		$1.53\!\!\times\!\!10^{-5}$ &  & $6.99\!\!\times\!\!10^{-15}$ &  \\
		$7.13\!\!\times\!\!10^{-5}$ &  & $2.85\!\!\times\!\!10^{-13}$ &  \\
		$1.92\!\!\times\!\!10^{-3}$ &  & $4.05\!\!\times\!\!10^{-13}$ &  \\
		\hline
		$5.25\!\!\times\!\!10^{3}$ & $3.04\!\!\times\!\!10^{-6}$ & $5.25\!\!\times\!\!10^{3}$ & $3.04\!\!\times\!\!10^{-6}$ \\
		& $6.29\!\!\times\!\!10^{-5}$ &  & $6.29\!\!\times\!\!10^{-5}$ \\
		& $3.51\!\!\times\!\!10^{-4}$ &  & $3.51\!\!\times\!\!10^{-4}$ \\
		& $6.26\!\!\times\!\!10^{-3}$ &  & $6.27\!\!\times\!\!10^{-3}$ \\
		& $9.53\!\!\times\!\!10^{-3}$ &  & $9.53\!\!\times\!\!10^{-3}$ \\
		& $4.42\!\!\times\!\!10^{0}$ &  & $4.42\!\!\times\!\!10^{0}$ \\
		& $4.59\!\!\times\!\!10^{0}$ &  & $4.59\!\!\times\!\!10^{0}$ \\
	\end{tabular}
\end{center}

\paragraph{Relaxation order $\omega=8$}
\begin{center}
	\begin{tabular}{ll|ll}
		\multicolumn{2}{c|}{$\varepsilon_{_{\mathrm{DIMACS}}}=10^{-6}$} &
		\multicolumn{2}{c}{$\varepsilon_{_{\mathrm{DIMACS}}}=10^{-16}$} \\
		\hline
		$X$ & $S$ & $X$ & $S$ \\
		\hline
		$2.11\!\!\times\!\!10^{-9}$ & $4.28\!\!\times\!\!10^{-13}$ & $2.80\!\!\times\!\!10^{-12}$ & $2.52\!\!\times\!\!10^{-16}$ \\
		$2.47\!\!\times\!\!10^{-9}$ &  & $3.99\!\!\times\!\!10^{-13}$ &  \\
		$7.53\!\!\times\!\!10^{-7}$ &  & $1.22\!\!\times\!\!10^{-13}$ &  \\
		$1.80\!\!\times\!\!10^{-6}$ &  & $7.28\!\!\times\!\!10^{-17}$ &  \\
		$1.67\!\!\times\!\!10^{-5}$ &  & $7.86\!\!\times\!\!10^{-16}$ &  \\
		$1.33\!\!\times\!\!10^{-4}$ &  & $9.71\!\!\times\!\!10^{-15}$ &  \\
		$1.19\!\!\times\!\!10^{-3}$ &  & $6.64\!\!\times\!\!10^{-13}$ &  \\
		$1.94\!\!\times\!\!10^{-2}$ &  & $1.82\!\!\times\!\!10^{-12}$ &  \\
		\hline
		$2.98\!\!\times\!\!10^{4}$ & $3.33\!\!\times\!\!10^{-7}$ & $2.98\!\!\times\!\!10^{4}$ & $3.34\!\!\times\!\!10^{-7}$ \\
		& $7.59\!\!\times\!\!10^{-6}$ &  & $7.60\!\!\times\!\!10^{-6}$ \\
		& $4.76\!\!\times\!\!10^{-5}$ &  & $4.77\!\!\times\!\!10^{-5}$ \\
		& $3.94\!\!\times\!\!10^{-4}$ &  & $3.95\!\!\times\!\!10^{-4}$ \\
		& $4.21\!\!\times\!\!10^{-3}$ &  & $4.21\!\!\times\!\!10^{-3}$ \\
		& $9.48\!\!\times\!\!10^{-3}$ &  & $9.49\!\!\times\!\!10^{-3}$ \\
		& $4.41\!\!\times\!\!10^{0}$ &  & $4.41\!\!\times\!\!10^{0}$ \\
		& $5.52\!\!\times\!\!10^{0}$ &  & $5.52\!\!\times\!\!10^{0}$ \\
	\end{tabular}
\end{center}

\paragraph{Relaxation order $\omega=9$}
\begin{center}
	\begin{tabular}{ll|ll}
		\multicolumn{2}{c|}{$\varepsilon_{_{\mathrm{DIMACS}}}=10^{-6}$} &
		\multicolumn{2}{c}{$\varepsilon_{_{\mathrm{DIMACS}}}=10^{-16}$} \\
		\hline
		$X$ & $S$ & $X$ & $S$ \\
		\hline
		$6.56\!\!\times\!\!10^{-10}$ & $5.41\!\!\times\!\!10^{-14}$ &
		$4.64\!\!\times\!\!10^{-20}$ & $1.54\!\!\times\!\!10^{-24}$ \\
		$1.22\!\!\times\!\!10^{-9}$  & & $1.08\!\!\times\!\!10^{-19}$ & \\
		$3.23\!\!\times\!\!10^{-7}$  & & $6.04\!\!\times\!\!10^{-17}$ & \\
		$6.31\!\!\times\!\!10^{-7}$  & & $9.94\!\!\times\!\!10^{-17}$ & \\
		$9.07\!\!\times\!\!10^{-6}$  & & $1.56\!\!\times\!\!10^{-15}$ & \\
		$4.53\!\!\times\!\!10^{-5}$  & & $5.80\!\!\times\!\!10^{-15}$ & \\
		$3.76\!\!\times\!\!10^{-4}$  & & $6.76\!\!\times\!\!10^{-14}$ & \\
		$2.69\!\!\times\!\!10^{-3}$  & & $4.45\!\!\times\!\!10^{-13}$ & \\
		$6.96\!\!\times\!\!10^{-2}$  & & $1.15\!\!\times\!\!10^{-11}$ & \\
		\hline
		$1.69\!\!\times\!\!10^{5}$   & $3.92\!\!\times\!\!10^{-8}$ &
		$1.69\!\!\times\!\!10^{5}$   & $3.93\!\!\times\!\!10^{-8}$ \\
		& $9.82\!\!\times\!\!10^{-7}$ & & $9.84\!\!\times\!\!10^{-7}$ \\
		& $6.19\!\!\times\!\!10^{-6}$ & & $6.20\!\!\times\!\!10^{-6}$ \\
		& $6.70\!\!\times\!\!10^{-5}$ & & $6.71\!\!\times\!\!10^{-5}$ \\
		& $2.77\!\!\times\!\!10^{-4}$ & & $2.78\!\!\times\!\!10^{-4}$ \\
		& $4.92\!\!\times\!\!10^{-3}$ & & $4.93\!\!\times\!\!10^{-3}$ \\
		& $6.54\!\!\times\!\!10^{-3}$ & & $6.55\!\!\times\!\!10^{-3}$ \\
		& $5.38$              & & $5.38$ \\
		& $5.49$              & & $5.49$ \\
	\end{tabular}
\end{center}

\subsection*{Problem~2}


\begin{center}
	{\small
		\begin{tabular}{ccccccc}
			\toprule
			$\omega$ & $m$ & matrix & $\varepsilon_{_{\mathrm{DIMACS}}}$ & rank$_+$ & $[\lambda_{\min+},\,\lambda_{\max+}]$ & $\kappa_+$ \\
			\midrule
			\multirow{2}{*}{6} & \multirow{2}{*}{28} & \multirow{2}{*}{$X$} & $10^{-6}$  & 1 & $[1.26\!\!\times\!\!10^{2},\,1.26\!\!\times\!\!10^{2}]$ & $1$ \\
			& & & $10^{-16}$ & 1 & $[1.26\!\!\times\!\!10^{2},\,1.26\!\!\times\!\!10^{2}]$ & $1$ \\
			\cmidrule(l){3-7}
			& & \multirow{2}{*}{$S$} & $10^{-6}$  & 27 & $[2.02\!\!\times\!\!10^{-7},\,7.51\!\!\times\!\!10^{0}]$ & $3.72\!\!\times\!\!10^{7}$ \\
			& & & $10^{-16}$ & 27 & $[2.02\!\!\times\!\!10^{-7},\,7.51\!\!\times\!\!10^{0}]$ & $3.72\!\!\times\!\!10^{7}$ \\
			\midrule
			\multirow{2}{*}{7} & \multirow{2}{*}{36} & \multirow{2}{*}{$X$} & $10^{-6}$  & 36 & $[4.77\!\!\times\!\!10^{-10},\,4.00\!\!\times\!\!10^{2}]$ & $8.38\!\!\times\!\!10^{11}$ \\
			& & & $10^{-16}$ & 5 & $[4.52\!\!\times\!\!10^{0},\,4.03\!\!\times\!\!10^{2}]$ & $8.92\!\!\times\!\!10^{1}$ \\
			\cmidrule(l){3-7}
			& & \multirow{2}{*}{$S$} & $10^{-6}$  & 36 & $[3.61\!\!\times\!\!10^{-11},\,7.78\!\!\times\!\!10^{0}]$ & $2.15\!\!\times\!\!10^{11}$ \\
			& & & $10^{-16}$ & 31 & $[7.56\!\!\times\!\!10^{-8},\,7.77\!\!\times\!\!10^{0}]$ & $1.03\!\!\times\!\!10^{8}$ \\
			\midrule
			\multirow{2}{*}{8} & \multirow{2}{*}{45} & \multirow{2}{*}{$X$} & $10^{-6}$  & 37 & $[4.85\!\!\times\!\!10^{-9},\,3.29\!\!\times\!\!10^{3}]$ & $6.77\!\!\times\!\!10^{11}$ \\
			& & & $10^{-16}$ & 7 & $[6.27\!\!\times\!\!10^{-1},\,3.30\!\!\times\!\!10^{3}]$ & $5.26\!\!\times\!\!10^{3}$ \\
			\cmidrule(l){3-7}
			& & \multirow{2}{*}{$S$} & $10^{-6}$  & 44 & $[5.85\!\!\times\!\!10^{-11},\,1.16\!\!\times\!\!10^{1}]$ & $1.98\!\!\times\!\!10^{11}$ \\
			& & & $10^{-16}$ & 38 & $[1.41\!\!\times\!\!10^{-9},\,1.16\!\!\times\!\!10^{1}]$ & $8.26\!\!\times\!\!10^{9}$ \\
			\midrule
			\multirow{2}{*}{9} & \multirow{2}{*}{55} & \multirow{2}{*}{$X$} & $10^{-6}$  & 51 & $[1.77\!\!\times\!\!10^{-8},\,9.84\!\!\times\!\!10^{3}]$ & $5.55\!\!\times\!\!10^{11}$ \\
			& & & $10^{-16}$ & 12 & $[6.91\!\!\times\!\!10^{-1},\,9.89\!\!\times\!\!10^{3}]$ & $1.43\!\!\times\!\!10^{4}$ \\
			\cmidrule(l){3-7}
			& & \multirow{2}{*}{$S$} & $10^{-6}$  & 48 & $[1.48\!\!\times\!\!10^{-11},\,1.19\!\!\times\!\!10^{1}]$ & $8.08\!\!\times\!\!10^{11}$ \\
			& & & $10^{-16}$ & 43 & $[7.70\!\!\times\!\!10^{-10},\,1.19\!\!\times\!\!10^{1}]$ & $1.55\!\!\times\!\!10^{10}$ \\
			\bottomrule
	\end{tabular}}
\end{center}

\subsection*{Problem~3}

\begin{center}
	{\small
		\begin{tabular}{ccccccc}
			\toprule
			$\omega$ & $m$ & matrix & $\varepsilon_{_{\mathrm{DIMACS}}}$ & rank$_+$ & $[\lambda_{\min+},\,\lambda_{\max+}]$ & $\kappa_+$ \\
			\midrule
			\multirow{2}{*}{4} & \multirow{2}{*}{15} & \multirow{2}{*}{$X$} & $10^{-6}$  & 15 & $[6.58\!\!\times\!\!10^{-9},\,2.40\!\!\times\!\!10^{3}]$ & $3.65\!\!\times\!\!10^{11}$ \\
			& & & $10^{-16}$ & 11 & $[1.48\!\!\times\!\!10^{-1},\,2.39\!\!\times\!\!10^{3}]$ & $1.62\!\!\times\!\!10^{4}$ \\
			\cmidrule(l){3-7}
			& & \multirow{2}{*}{$S$} & $10^{-6}$  & 12 & $[1.96\!\!\times\!\!10^{-12},\,1.00\!\!\times\!\!10^{0}]$ & $5.10\!\!\times\!\!10^{11}$ \\
			& & & $10^{-16}$ & 4 & $[3.01\!\!\times\!\!10^{-5},\,1.00\!\!\times\!\!10^{0}]$ & $3.33\!\!\times\!\!10^{4}$ \\
			\midrule
			\multirow{2}{*}{5} & \multirow{2}{*}{21} & \multirow{2}{*}{$X$} & $10^{-6}$  & 20 & $[1.01\!\!\times\!\!10^{-6},\,8.11\!\!\times\!\!10^{4}]$ & $7.99\!\!\times\!\!10^{10}$ \\
			& & & $10^{-16}$ & 16 & $[7.98\!\!\times\!\!10^{-2},\,8.17\!\!\times\!\!10^{4}]$ & $1.02\!\!\times\!\!10^{6}$ \\
			\cmidrule(l){3-7}
			& & \multirow{2}{*}{$S$} & $10^{-6}$  & 17 & $[1.48\!\!\times\!\!10^{-12},\,1.00\!\!\times\!\!10^{0}]$ & $6.74\!\!\times\!\!10^{11}$ \\
			& & & $10^{-16}$ & 5 & $[5.00\!\!\times\!\!10^{-7},\,1.00\!\!\times\!\!10^{0}]$ & $2.00\!\!\times\!\!10^{6}$ \\
			\midrule
			\multirow{2}{*}{6} & \multirow{2}{*}{28} & \multirow{2}{*}{$X$} & $10^{-6}$  & 26 & $[6.93\!\!\times\!\!10^{-6},\,2.71\!\!\times\!\!10^{6}]$ & $3.91\!\!\times\!\!10^{11}$ \\
			& & & $10^{-16}$ & 22 & $[5.02\!\!\times\!\!10^{-2},\,2.79\!\!\times\!\!10^{6}]$ & $5.55\!\!\times\!\!10^{7}$ \\
			\cmidrule(l){3-7}
			& & \multirow{2}{*}{$S$} & $10^{-6}$  & 18 & $[1.65\!\!\times\!\!10^{-12},\,1.00\!\!\times\!\!10^{0}]$ & $6.06\!\!\times\!\!10^{11}$ \\
			& & & $10^{-16}$ & 6 & $[9.67\!\!\times\!\!10^{-9},\,1.00\!\!\times\!\!10^{0}]$ & $1.03\!\!\times\!\!10^{8}$ \\
			\bottomrule
	\end{tabular}}
\end{center}

%


\end{appendices}

\section*{Acknowledgements}
The author acknowledges the financial support of the Czech Ministry of Education, Youth and Sports project ROBOPROX (reg.\ no.\ CZ.02.01.01/00/22\_008/0004590) that was co-funded by the European Union.

The author would like to thank Lorenzo Baldi and Didier Henrion for stimulating this work and for providing the optimization problems considered in this paper.

The author used AI to discuss the mathematical exposition, flag possible logical inconsistencies, and help revise and clarify the manuscript.
The author independently evaluated all such suggestions and verified
the mathematical arguments, numerical data, and conclusions. The
author remains fully responsible for the content of the manuscript.

\bibliography{dpop}

\end{document}